\documentclass[12pt]{article}

\usepackage{arxiv}
\usepackage{authblk}
\usepackage{setspace}
\usepackage[utf8]{inputenc} 
\usepackage[T1]{fontenc}    
\usepackage{hyperref}       
\usepackage{url}            
\usepackage{booktabs}       
\usepackage{amsfonts}       
\usepackage{nicefrac}       
\usepackage{microtype}      
\usepackage{lipsum}
\usepackage{graphicx}
\usepackage{graphics}
\usepackage{amsmath}
\usepackage{caption}
\usepackage{subcaption}
\usepackage{amssymb}
\usepackage{algorithm}
\usepackage{tcolorbox} 
\usepackage{mathtools}
\usepackage{booktabs}
\usepackage{mathrsfs}
\usepackage{enumerate}
\usepackage{epstopdf}
\usepackage{yhmath}
\usepackage[named]{algo}
\usepackage[titletoc]{appendix}
\usepackage{color}
\usepackage[numbers,sort&compress]{natbib}
\usepackage{hyperref}
\usepackage{verbatim}
\graphicspath{ {./images/} }

\newtheorem{lemma}{Lemma}

\newtheorem{remark}{Remark}
\newtheorem{theorem}{Theorem}
\newtheorem{definition}{Definition}

\newtheorem{proposition}{Proposition}
\title{Elongation of Curvature-Bounded Path}

\author[1,2*]{Zheng Chen}
\author[1]{Kun Wang}
\author[1]{Yang Lu}
\affil[1]{School of Aeronautics and Astronautics, Zhejiang University, Hangzhou 310027, Zhejiang, China}
\affil[2]{State Key Laboratory of Fluid Power and Mechatronic Systems, Hangzhou 310027, Zhejiang, China}
\affil[*]{Email address: z-chen@zju.edu.cn}

\begin{document}
\maketitle
\begin{abstract}
The paper is concerned with elongating the shortest curvature-bounded path  between two oriented points to an expected length. The elongation of  curvature-bounded paths to an expected length is fundamentally important to plan missions for nonholonomic-constrained vehicles in many practical applications, such as coordinating multiple nonholonomic-constrained vehicles to reach a destination simultaneously or performing a mission with a strict time window. In the paper, 
 the explicit conditions for the existence of curvature-bounded paths joining  two oriented points with an expected length are established by applying the properties of the reachability set of curvature-bounded paths. These existence conditions are numerically verifiable, allowing readily checking the existence of curvature-bounded paths between two prescribed oriented points with a desired length. In addition, once the existence conditions are met, elongation strategies are provided in the paper to get curvature-bounded paths with expected lengths. Finally, some examples of minimum-time path planning for multiple fixed-wing aerial vehicles to cooperatively achieve a triangle-shaped flight formation are presented, illustrating and verifying the developments of the paper.

\end{abstract}

\keywords{Dubins vehicle\and Curvature-bounded path\and Cooperative guidance\and Path elongation}

\section{Introduction}

The model of unidirectional nonholonomic vehicles, moving at a constant speed with a minimum turn radius, provides a very good approximation to the kinematics of a large class of vehicles, such as  fixed-wing aerial vehicles, autonomous underwater vehicles, and uninhabited ground vehicles, just to name a few. Following the work by A. A. Markov \cite{markov1887some} in 1887 and the work  by L. E. Dubins \cite{Dubins:1957} in 1957, such a nonholonomic vehicle has been named Markov-Dubins vehicle (see, e.g., \cite{kaya2017markov,bakolas2013optimal}) or simply Dubins vehicle (see, e.g., \cite{Bui:1994,meyer2015dubins,boissonnat1994shortest}). Since the Dubins vehicle has a minimum turn radius, it follows that the curvature of any feasible path of Dubins vehicle is bounded almost everywhere. For the sake of simplicity, we use the term of 
"curvature-bounded path" to denote the feasible path of Dubins vehicle in this paper.

Up to present, the minimum-time curvature-bounded paths with various constraints have been extensively studied in the literature. Note that the minimum-time path is the same as the shortest path as the speed is constant. In the seminal paper \cite{Dubins:1957} by L. E. Dubins, the shortest curvature-bounded path between two oriented points in the tangent bundle of $\mathbb{R}^2$ (each oriented point consists of a $2$-dimenstional point and a tangent vector at the point) was studied by geometric analysis, showing that the shortest curvature-bounded path between two oriented points can be computed within a constant time by comparing at most 6 candidate paths. With the advent of optimal control theory, this result was proven in an alternative way by combining Pontryagin's maximum principle \cite{pontryagin1987mathematical} and geometric techniques \cite{boissonnat1994shortest,sussmann1991shortest}. Considering to start from an oriented point and to reach a point instead of an oriented point (i.e., the direction of the final tangent vector is not constrained), it was proven in \cite{Bui:1994} that the shortest curvature-bounded path  must lie in a sufficient family of 4 candidate paths.

Due to the importance of using curvature-bounded paths in real-world scenarios, the shortest curvature-bounded path with more complex environmental and boundary constraints have been widely studied. Examples include, but not limited to, shortest curvature-bounded paths passing through multiple waypoints \cite{chen2019shortest,ma2006receding,chen2019relaxed}, passing through multiple regions \cite{chen2021descent,vavna2015dubins,vavna2018dubins}, encircling a target  \cite{chen2020dubins,manyam2019optimal}, avoiding obstacles \cite{jha2020shortest,boissonnat1996polynomial,giordano2009shortest,savkin_hoy_2013},   intercepting a moving target \cite{zheng2021time,zheng2021time_angle,Buzikov:2021,meyer2015dubins}, and moving in tunnel-like environments \cite{MATVEEV2020108831} or in uniform current drift \cite{bakolas2013optimal,techy2009minimum,mcgee2007optimal}. 

In addition to the above mentioned shortest curvature-bounded paths with various environmental and boundary constraints, it is also quite important to find curvature-bounded paths with expected lengths in practical applications \cite{ding2019curvature}. For instance, coordinating multiple nonholonomic vehicles from different initial conditions to the same final condition simultaenously requires finding   curvature-bounded paths with an expected length \cite{Yao:2017,Yao:2020Trajectory,SHANMUGAVEL20101084}. Another example is that a Dubins vehicle performes a mission with a strict time window, which requires planning a curvature-bounded path with a specific length. 
For this reasion, the issue of finding curvature-bounded paths with   expected lengths has attracted extensive attention in the past decades. 

In order to find a curvature-bounded path with an expected length, it is common in the literture to first find the shortest curvature-bounded path with the same boundary conditions, and then to use some elongation strategies to elongate the shortest curvature-bounded path to the expected length.  When the direction of the final tangent vector is not constrained, by dividing the $2$-dimensional plane into $5$ subregions, some iterative algorithms were developed in \cite{schumacher2003path} to elongate curvature-bounded paths depending on the locations of the final point in the five subregions. Meyer, Isaiah, and Shima   \cite{meyer2015dubins} proposed three strategies to elongate curvature-bounded paths for intercepting a moving target at a given time. Recently, it was proven by Ding, Xin, and Chen \cite{ding2019curvature} that if the boundary conditions lie in a specific set, the shortest curvature-bounded path  cannot be elongated to arbitrary length. 

In all the papers mentioned in the previous paragraph, it is assumed that the direction of the final tangent vector is free. When the tangent vectors at both endpoints are fixed, the shortest curvature-bounded paths are composed by circular arcs and straight line segments \cite{Dubins:1957,boissonnat1994shortest}; to be specific, the geometric pattern of each shortest curvature-bounded path between two oriented points is CCC, CSC, or their substrings, where ``C'' denotes a circular arc and ``S'' denotes a straight line segment. It was proposed in \cite{shanmugavel2005path} to elongate the shortest curvature-bounded path between two oriented points by increasing the turning radius of circular arcs  and a bisection method was used to find the optimal turning radii of paths.  
Ortiz, Kingston, and Langbort  \cite{ortiz2013multi} provided path elongation strategies with a strict assumption that the shortest curvature-bounded paths take a geometric pattern of CSC. It is worth mentioning that the Dubins paths with clothoid arcs were used in \cite{SHANMUGAVEL20101084} to generate curvature-bounded paths with expected lengths.

Considering that the final point lies in a manifold, a homotopy method was used to generate curvature-bounded paths with expected lengths \cite{Yao:2017}. In addition, it was shown in  \cite{Yao:2017} that, given two specific oriented points, there will be a non-zero interval determined by the two oriented points, so that the length of any curvature-bounded path between the two oriented points does not lie in the interval; this means that the shortest curvature-bounded path between the two oriented points cannot be elongated to arbitrary length. Recently, without  requiring the length of the curvature-bounded path to be strictly equal to an expected value, a homotopy method was used in \cite{Yao:2020} to elongate a curvature-bounded path to a length as close  as possible to the expected length.

Topologic techniques were employed by J. Ayala {\it et al.} in \cite{Ayala+2017+283+292,Ayala:2016,Ayala:2015,Ayala:2014} to show that the curvature-bounded paths between two prescribed oriented points lie in two homotopy classes. If the two oriented points take some specific values, the two homotopy classes do not connect with each other. This coincides with the results obtained by W. Yao {\it et al.} in \cite{Yao:2017,Yao:2020} that the shortest curvature-bounded path cannot be elongated to arbitrary length. All in all, it is clear according to the cited papers \cite{Yao:2017,Yao:2020,Ayala+2017+283+292,Ayala:2016,Ayala:2015,Ayala:2014} that given two oriented points and an expected length, one may not be able to find curvature-bounded paths between the two oriented points with the expected length.

However, to the authors' best knowledge, given two oriented points, it is not clear in the literature how to compute the exact intervals so that for every expected length in the intervals there exits a curvature-bounded path joining the two oriented points. Finding 
such intervals is a fundamental issue for evaluating the existence of curvature-bounded paths between two oriented points with an expected length. In this paper,  the properties of the reachability set constructed by Patsko, Pyatko, and Fedotov in the remarkable work \cite{patsko2003three} are employed to 
establish the exact intervals, which will further give rise to the explicit conditions for the existence of curvature-bounded path with an expected length. 

For the completeness of this paper, it is first proven  that if the shortest curvature-bounded path between two oriented points takes a geometric pattern of CCC, it can be enlongated to arbitrary length (cf. Proposition  \ref{LE:CCC}). In the case that the shortest curvature-bounded path is of CSC, we devide the boundary conditions into two seperate sets $\mathcal{O}$ and $\nabla \mathcal{O}$. If the two endpoints lie in $\mathcal{O}$, the shortest curvature-bounded path can be elongated to arbitrary length (cf. Proposition \ref{LE:first_set}). If the two endpoints lie in  $\nabla \mathcal{O}$, there exists an interval so that the length of any curvature-bounded path between the two endpoints is not in the interval (cf. Proposition \ref{LE:Lm_L1} and Proposition \ref{LE:non-existence}). The boundary of the interval is explicitly devised. As a result, given any two oriented points and any expected length, one can readily check the existence of a curvature-bounded path joining the two oriented points with the expected length. Once it exists, an elongation strategy is provided in the paper to elongate the shortest curvature-bounded path to  the expected length. 

The paper is organized as follows. Definitions and notations are presented in Section \ref{SE:Preliminary}. Given any expected length, the conditions for the existence of a curvature-bounded path  are established in Section \ref{SE:Elongation}, and once the existence conditions are met, an elongation strategy is provided. Section \ref{SE:Numerical} presents some  numerical examples to illustrate the developments of the paper, and this paper concludes by Section \ref{SE:Conclusions}.
\section{Definitions and Notations}\label{SE:Preliminary}
Denote by $T\mathbb{R}^2$ the tangent bundle of $\mathbb{R}^2$. Any element in $T\mathbb{R}^2 $ corresponds to an oriented point $(\boldsymbol{x},\boldsymbol{v})$ where $\boldsymbol{x}$ is a point in $\mathbb{R}^2$ and $\boldsymbol{v}$ is a tangent vector to $\mathbb{R}^2$ at $\boldsymbol{x}$.  For notational simplicity, we also use capital letters $X$ and $Y$ to denote the elements in $T\mathbb{R}^2$ in this paper, and we set 
$$X \coloneqq (\boldsymbol{x},\boldsymbol{v})\ \text{and}\ Y \coloneqq (\boldsymbol{y},\boldsymbol{w}).$$

\begin{definition}[Curvature-Bounded Path \cite{Ayala:2016}]
Given any two oriented points  $(\boldsymbol{x},\boldsymbol{v})$ and $(\boldsymbol{y},\boldsymbol{w})$ in $T\mathbb{R}^2$, a path $\gamma : [0,s]\rightarrow \mathbb{R}^2$ connecting the two oriented points is a curvature-bounded path if
\begin{itemize}
\item $\gamma$ is $C^1$ and piecewise $C^2$;
\item $\gamma$ is parameterised by arc length, i.e., $\|\gamma^{\prime}(t)\| = 1$ for all $t\in [0,s]$;
\item $\gamma(0) = \boldsymbol{x}$, $\gamma^{\prime}(0) = \boldsymbol{v}$, $\gamma(s) = \boldsymbol{y}$, and $\gamma^{\prime}(s) = \boldsymbol{w}$;
\item $\gamma^{\prime\prime}(t)\leq \kappa$ for all $t\in [0,s]$ when defined, where $\kappa> 0$ is a constant. 
\end{itemize}
\end{definition}
\begin{definition}
Given any two oriented points $X$ and $Y$ in $T\mathbb{R}^2$  such that $X\neq Y$, we denote by $\Gamma(X,Y)$ the space of all the curvature-bounded paths from $X$ to $Y$. 
\end{definition}
\begin{definition}
Given any two oriented points $X$ and $Y$ in $T\mathbb{R}^2$  such that $X\neq Y$,  for each path $\gamma \in \Gamma(X,Y)$,  we denote by $\ell(\gamma)$ the   length of $\gamma$.
\end{definition}
\begin{definition}\label{DE:Sm}
Given any two oriented points $X$ and $Y$ in $T\mathbb{R}^2$  such that $X\neq Y$, we denote by $\gamma_m$ the shortest curvature-bounded path from  $X$ to $Y$, and denote by $\ell_m>0$ the arc length of $\gamma_m$, i.e., $\ell_m = \ell(\gamma_m)$.
\end{definition}
\begin{definition}\label{DE:tangent}
Given any two oriented points $X$ and $Y$ in $T\mathbb{R}^2$  such that $X\neq Y$, a  path $\gamma \in \Gamma(X,Y)$ is said to have parallel tangents if there are two different points on $\gamma$ so that the tangent vectors at the two points are parallel with opposite directions. 
\end{definition}
It has been shown by L. E. Dubins  \cite{Dubins:1957} that the shortest curvature-bounded path between two oriented points lies in a sufficiently family of six candidates.  Denote by ``C'' a circular arc with radius of $1/\kappa$,  and denote by ``S'' a straight line segment. Then,  we have the following remark. 
\begin{remark}[Dubins path \cite{Dubins:1957}]\label{RE:Dubins}
Given any two oriented points $X$ and $Y$ in $T\mathbb{R}^2$ such that $X\neq Y$,  the shortest curvature-bounded path from $X$ to $Y$ takes a geometric pattern of  either $\mathrm{CCC}$ or $\mathrm{CSC}$ or their substrings where
\begin{itemize}
\item $\mathrm{CCC}$ = $\{$ $\mathrm{RLR}$, $\mathrm{LRL}$ $\}$,
\item $\mathrm{CSC}$ = $\{$ $\mathrm{RSR}$, $\mathrm{RSL}$, $\mathrm{LSR}$, $\mathrm{LSL}$$\}$
\end{itemize}
where $\mathrm{R}$ (resp. $\mathrm{L}$) means that the corresponding circular arc has a right-turning (resp. left-turning) direction.
\end{remark}
Thanks to the geometric patterns in Remark \ref{RE:Dubins}, the shortest curvature-bounded path can be analytically computed by comparing at most $6$ candidate paths; see \cite{Shkel:2001}. 
 
\begin{figure}[!htp]
\centering
\begin{subfigure}{4cm}
\includegraphics[width = 4cm]{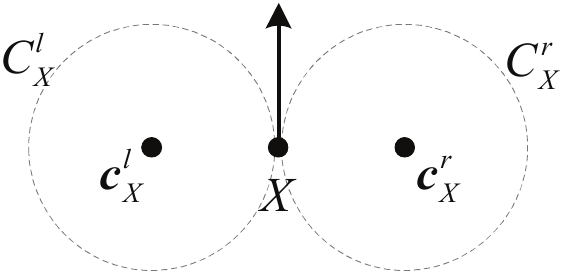}
\caption{$\boldsymbol{c}_{X}^r$ and $\boldsymbol{c}_X^l$}
\label{Fig:Center_X}
\end{subfigure}
\begin{subfigure}{4cm}
\includegraphics[width = 4cm]{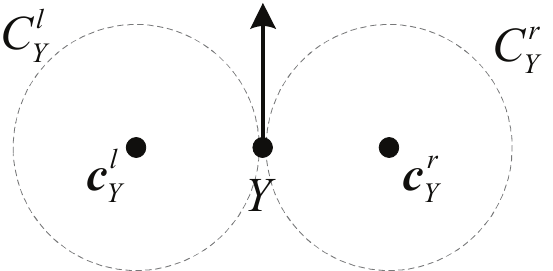}
\caption{$\boldsymbol{c}_{Y}^r$ and $\boldsymbol{c}_Y^l$}
\label{Fig:Center_Y}
\end{subfigure}
\caption{The geometry of $C_{X}^r$, $C_X^l$, $C_{Y}^r$, and $C_Y^l$ with respect to $X$ and $Y$.}
\end{figure}

Before proceeding, we present some useful notations. Let $C_X^r$ and $C_X^l$ be two tangent circles of radius $1/\kappa$,  lying on the right and left side of the initial oriented points $X$,  respectively, as shown by the two dashed circles in Fig.~\ref{Fig:Center_X}. We denote by $\boldsymbol{c}_X^r$ and $\boldsymbol{c}_X^l$ the centers of $C_X^r$ and $C_X^l$, respectively. The same applies to $C_Y^r$, $C_Y^l$, $\boldsymbol{c}_Y^r$, and $\boldsymbol{c}_Y^l$, as shown in Fig.~\ref{Fig:Center_Y}. Let $\mathrm{C}_{\eta}$ be a circular arc with radian of $\eta\geq 0$, and denote by $\mathrm{S}_d$ a straight line segment with length of $d\geq 0$.  With these new notations, when necessary we will represent CCC and CSC by C$_{\eta}$C$_{\zeta}$C$_{\xi}$ and C$_{\eta}$S$_{d}$C$_{\xi}$, respectively. 
\section{Elongation of Curvature-Bounded Path}\label{SE:Elongation}
The elongation of curvature-bounded path  is closely related to the existence of curvature-bounded paths with an expected length.  In this section, the conditions for the existence of curvature-bounded path with an expected length will be established, which will further give rise to the elongationability of curvature-bounded paths.
\subsection{Elongation of CCC-Path}
\begin{figure}[!htbp]
\centering
\includegraphics[width = 4cm]{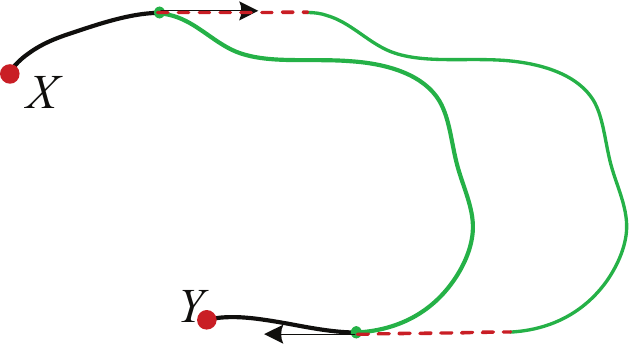}
\caption{Elongation of curvature-bounded path with parallel tangents.}
\label{Fig:Tangent}
\end{figure}

This subsection will show how to elongate the shortest curvature-bounded path $\gamma_m$ if its geometric pattern is of type CCC. Before proceeding, we first show that any curvature-bounded path with parallel tangents (cf. Definition \ref{DE:tangent}) can be elongated to arbitrary length by the following lemma. 
\begin{lemma}[J. Ayala \cite{Ayala+2017+283+292}]\label{LE:parallel}
Given any curvature-bounded path $\gamma \in \Gamma(X,Y)$, if it has parallel tangents, then for any $s\geq \ell(\gamma)$ there exists a curvature-bounded path $\bar{\gamma}\in \Gamma(X,Y)$ so that $s = \ell(\bar{\gamma})$.
\end{lemma}
The proof of this lemma can be found in \cite{Ayala+2017+283+292}, and an example for elongating the curvature-bounded path with parallel tangents is illustrated in Fig.~\ref{Fig:Tangent}. 
\begin{proposition}\label{LE:CCC}
Given any two oriented points $X$ and $Y$ in $T\mathbb{R}^2$  such that $X\neq Y$, if the shortest curvature-bounded path $\gamma_{m}$ takes a geometric pattern of  $ \mathrm{C_{\eta}C_{\zeta}C_{\xi}}$, then for any $s \geq \ell_m$ there exists  a curvature-bounded path $\gamma \in \Gamma(X,Y)$ such that $s = \ell(\gamma)$.
\end{proposition}
Proof. In view of \cite[Lemma 3]{Bui:1994}, if the shortest curvature-constrained path $\gamma_m$ is of type $\mathrm{C_{\eta}C_{\zeta}C_{\xi}}$, we have that   the middle circular arc is a major arc, i.e., $\zeta \in (\pi,2\pi)$. Note that a circular arc $\mathrm{C}_{\zeta}$ with $\zeta \geq \pi$ admits parallel tangents. Thus, according to Lemma \ref{LE:parallel}, for every $s\geq \ell_m$ there exists a curvature-bounded path $\gamma \in \Gamma(X,Y)$ so that $s = \ell(\gamma)$, completing the proof.$\square$\\

\begin{figure}[!htp]
\centering
\includegraphics[width = 4cm]{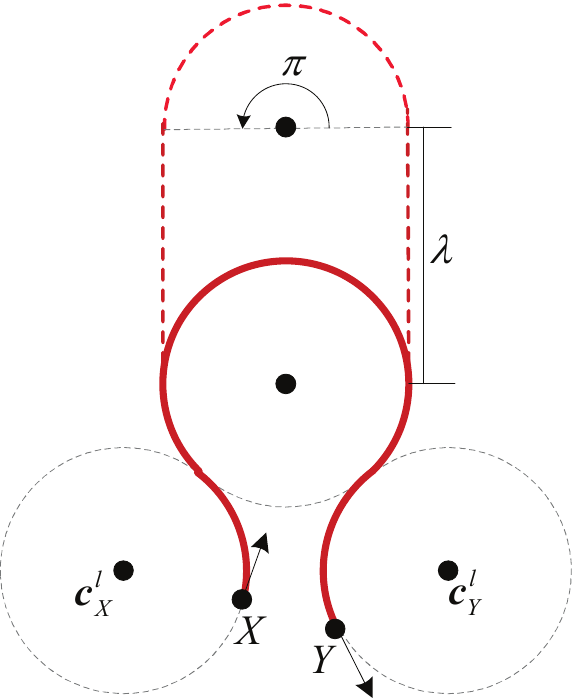}
\caption{An elongation strategy for the shortest path of type CCC.}
\label{Fig:CCC_Elongation}
\end{figure}

Proposition \ref{LE:CCC} indicates that if the shortest curvature-bounded path is of type CCC, it can be elongated to arbitrary length without breaking the constraint on culvature.  An elongation strategy for a CCC-path is illustrated in Fig.~\ref{Fig:CCC_Elongation} where the path is elongated by changing the value of $\lambda\in [0,+\infty)$. In the following subsection, we shall show the elongationability of the shortest curvature-bounded path if its geometric pattern is of type CSC. 
\subsection{Elongation of CSC-Path}
\begin{figure}[!htp]
\centering
\includegraphics[width = 4cm]{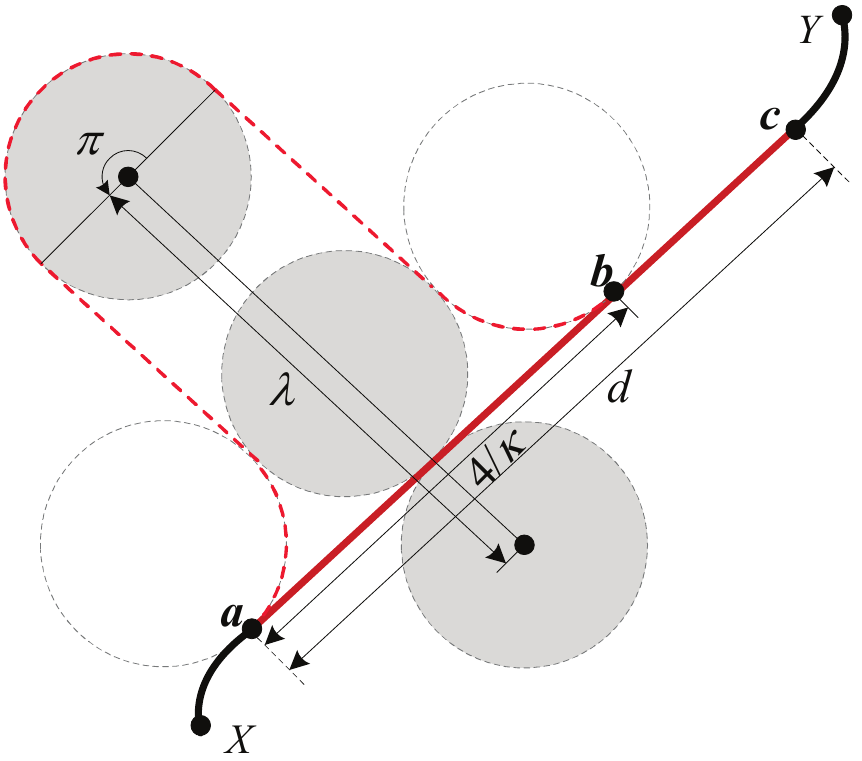}
\caption{Elongation of a path in $\Gamma(X,Y)$ with  a straight line segment and the length of the straight line segment is no less than  $4/\kappa$.}
\label{Fig:Straight}
\end{figure}

For notational simplicity, we define the following five sets. 
\begin{align}
\begin{split}
\mathcal{O}_1 =&\ \{(X,Y)\in (T\mathbb{R}^2)^2 | \gamma_m \in \mathrm{C}_{\eta}\mathrm{S}_{d}\mathrm{C}_{\xi}\ \text{with}\ \eta\geq \pi\}
\end{split}\nonumber\\
\begin{split}
\mathcal{O}_2 =&\ \{(X,Y)\in (T\mathbb{R}^2)^2 | \gamma_m \in \mathrm{C}_{\eta}\mathrm{S}_{d}\mathrm{C}_{\xi}\ \text{with}\ \xi\geq \pi\}
\end{split}\nonumber\\
\begin{split}
\mathcal{O}_3 =&\ \{(X,Y)\in (T\mathbb{R}^2)^2 | \gamma_m \in \mathrm{C}_{\eta}\mathrm{S}_{d}\mathrm{C}_{\xi}\ \text{with}\ d \geq \frac{4}{\kappa}\}
\end{split}\nonumber\\
\begin{split}
\mathcal{O}_4 =&\  \{(X,Y)\in (T\mathbb{R}^2)^2 | \gamma_m \in \mathrm{C}_{\eta}\mathrm{S}_{d}\mathrm{C}_{\xi}\ \text{with}\  d (\boldsymbol{c}_X^r,\boldsymbol{c}_Y^r) \geq 4/\kappa\}
\end{split}\nonumber\\
\begin{split}
\mathcal{O}_5 =&\  \{(X,Y)\in (T\mathbb{R}^2)^2 | \gamma_m \in \mathrm{C}_{\eta}\mathrm{S}_{d}\mathrm{C}_{\xi}\ \text{with}\  d (\boldsymbol{c}_X^l,\boldsymbol{c}_Y^l) \geq 4/\kappa\}
\end{split}\nonumber
\end{align}
where the function $d:\mathbb{R}^2\times\mathbb{R}^2\rightarrow \mathbb{R}_+$ denotes the Euclidean distance between two points. 
Let $\mathcal{O}\subset (T\mathbb{R}^2)^2$ be the union of $\mathcal{O}_1$, $\mathcal{O}_2$, $\ldots$, $\mathcal{O}_5$, i.e.,
$$\mathcal{O} \coloneqq \mathcal{O}_1\cup \mathcal{O}_2 \cup \mathcal{O}_3 \cup \mathcal{O}_4 \cup \mathcal{O}_5,$$ 
and let $\nabla \mathcal{O}\subset (T\mathbb{R}^2)^2$ be the complementary set of $\mathcal{O}$, i.e., 
\begin{align}
\begin{split}
\nabla{\mathcal{O}} \coloneqq & \{(X,Y)\in (T\mathbb{R}^2)^2| \gamma_m \in \mathrm{CSC}, \ (X,Y)\not \in \mathcal{O} \}.
\end{split}\nonumber
\end{align} 
In view of the definitions of $\mathcal{O}$ and $\nabla\mathcal{O}$, it is clear that for any $(X,Y)\in (T\mathbb{R}^2)^2$ so that the geometric pattern of $\gamma_m$ is CSC, we have $(X,Y)\in \mathcal{O}\cup \nabla\mathcal{O}$. By the following lemmas and theorems, we shall establish the  conditions for the existence of curvature-bounded path $\gamma \in \Gamma(X,Y)$ with expected lengths for $(X,Y)$ in the two seperate sets $\mathcal{O}$ and $\nabla \mathcal{O}$.
\begin{lemma}\label{LE:d4k}
Given any curvature-bounded path $\gamma \in \Gamma(X,Y)$, if it has a straight-line segment  and if the  length of the straight line segment is no less than $4/\kappa$, then for any $s\geq \ell(\gamma)$ there exists a curvature-bounded path $\bar{\gamma}\in \Gamma(X,Y)$ so that $s = \ell(\bar{\gamma})$.
\end{lemma}
Proof. Denote by $d\geq 4/\kappa$ the length of the straight line segment, as shown in Fig.~\ref{Fig:Straight}. Denote by $\boldsymbol{a}$ and $\boldsymbol{c}$ the initial and final points of the straight line segment.  Then, we can choose a point $\boldsymbol{b}$ on the straight line segment so that the distance between $\boldsymbol{a}$ and $\boldsymbol{b}$ is $4/\kappa$. We can use a circular disk of radius $1/\kappa$ to deform the straight line segment from $\boldsymbol{a}$ to $\boldsymbol{b}$ without changing the tangent vectors at $\boldsymbol{a}$ to $\boldsymbol{b}$, as shown by Fig.~\ref{Fig:Straight}. Thus, the path $\gamma$ can be elongated to arbitrary length without breaking the constraint on maximum curvature, completing the proof. $\Box$
\begin{proposition}\label{LE:first_set}
Given any two oriented points $X$ and $Y$ in $T\mathbb{R}^2$ so that $(X,Y)\in \mathcal{O}$,  for any $s\geq \ell_m$ there exists a curvature-bounded path $\gamma \in \Gamma(X,Y)$ such that $s = \ell(\gamma)$. 
\end{proposition}
Proof. According to the definitions of $ \mathcal{O}_1$ and $ \mathcal{O}_2$, as a major circular arc has parallel tangents, we have that the shortest curvature-bounded path $\gamma_m$ has parallel tangents if $(X,Y)\in \mathcal{O}_1\cup \mathcal{O}_2$. Therefore, according to Lemma \ref{LE:parallel}, the shortest curvature-bounded path $\gamma_m$ can be enlongated to arbitrary length. Therefore, if $(X,Y)\in \mathcal{O}_1\cup \mathcal{O}_2$, for any $s\geq \ell_m$ there exists a curvature-bounded path $\gamma \in \Gamma(X,Y)$ so that $s = \ell(\gamma)$. 

From now on, we consider $(X,Y)\in \mathcal{O}_3$. In this case, there exists a straight line segment with its length no less than $4/\kappa$ along   the shortest curvature-bounded path $\gamma_m$. Thus, according to  Lemma \ref{LE:d4k}, the shortest curvature-bounded path $\gamma_m$ can be elongated to arbitrary length, indicating that for any $s\geq \ell_m$ there exists a curvature-bounded path $\gamma \in \Gamma(X,Y)$ so that $s = \ell(\gamma)$. 

\begin{figure}[!htp]
\centering
\begin{subfigure}{3.9cm}
\centering
\includegraphics[width =3.7cm]{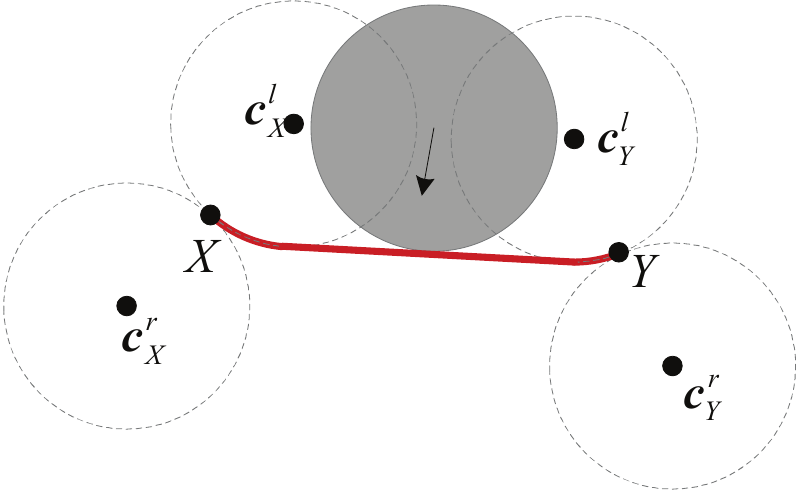}
\caption{}
\label{Fig:d4_Elongation1}
\end{subfigure}
\begin{subfigure}{3.9cm}
\centering
\includegraphics[width = 3.7cm]{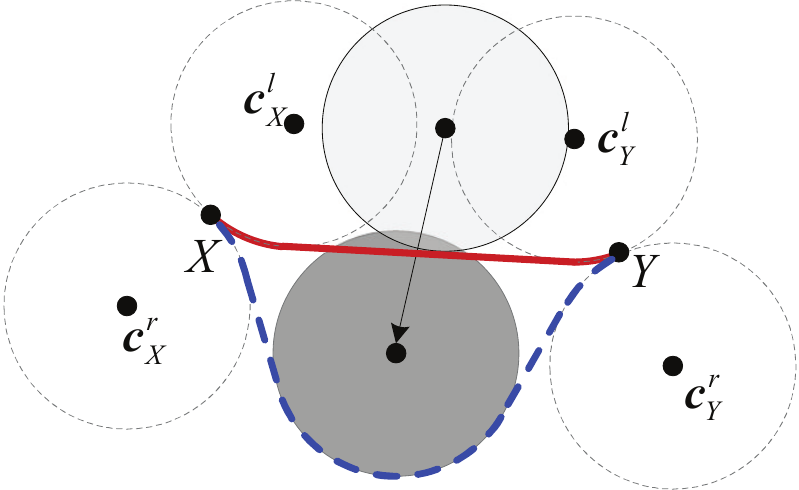}
\caption{}
\label{Fig:d4_Elongation2}
\end{subfigure}
\caption{An elongation strategy for $\gamma_m\in \Gamma(X,Y)$ with  $(X,Y)\in \mathcal{O}_4$.}
\label{Fig:d4_Elongation}
\end{figure}

Regarding the case that $(X,Y)\in \mathcal{O}_4$, as shown in Fig.~\ref{Fig:d4_Elongation1}, we are able to deform the shortest curvature-bounded path (red solid curve) by moving a circular disk of radius $1/\kappa$ without breaking the constraint on curvature, as shown by the blue dashed curve in Fig.~\ref{Fig:d4_Elongation2}. Thus, in the case of $(X,Y)\in \mathcal{O}_4$, the shortest curvature-bounded path $\gamma_m$ can be elongated to arbitrary length. The same strategy can be applied to elongate the shortest curvature-bounded path $\gamma_m$ if $(X,Y)\in \mathcal{O}_5$. Thus, if $(X,Y)\in \mathcal{O}_4\cup \mathcal{O}_5$, for any $s\geq \ell_m$ there exists a curvture-bounded path $\gamma\in \Gamma(X,Y)$ so that $s = \ell(\gamma)$, completing the proof. $\Box$

As a result of Proposition \ref{LE:first_set}, it is apparent that if $(X,Y)\in \mathcal{O}$, the shortest curvature-bounded path $\gamma_m$ can be enlongated to arbitrary length without breaking the constraint on maximum curvature. 

\begin{figure}[!htp]
\centering
\begin{subfigure}{5cm}
\includegraphics[width = 5cm]{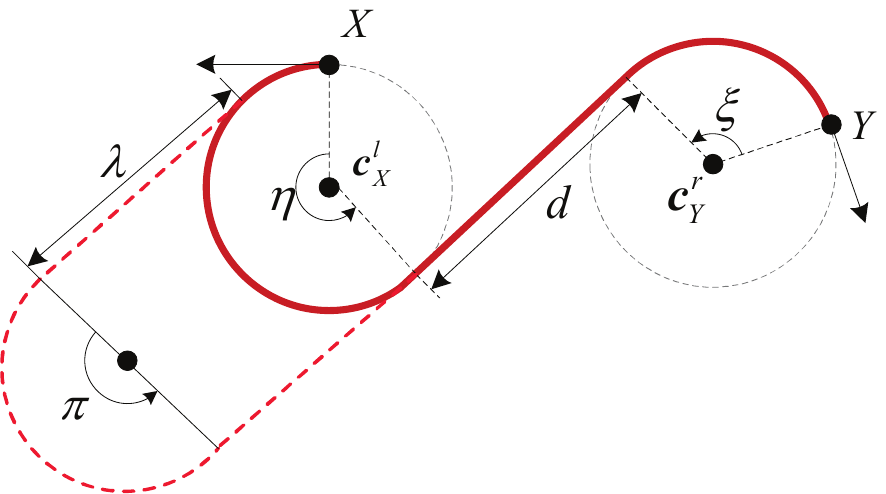}
\caption{$\mathrm{C_{\eta}S_d C_{\xi}}$ with $\eta > \pi$}
\label{Fig:LSR_Elongation}
\end{subfigure}
\begin{subfigure}{5cm}
\includegraphics[width = 5cm]{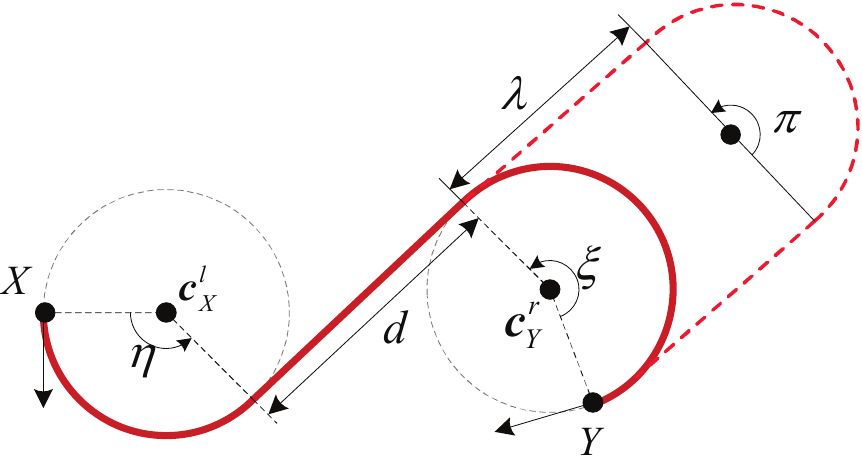}
\caption{$\mathrm{C_{\eta}S_d C_{\xi}}$ with $\xi > \pi$}
\label{Fig:LSR1_Elongation}
\end{subfigure}
\begin{subfigure}{5cm}
\includegraphics[width = 4.8cm]{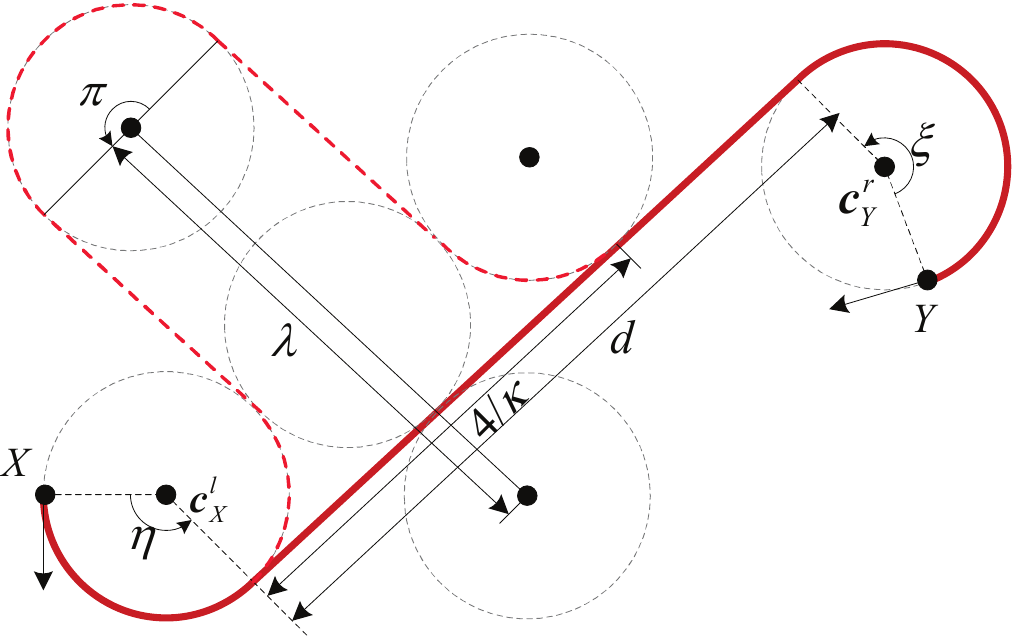}
\caption{$\mathrm{C_{\eta}S_d C_{\xi}}$ with $d > 4/\kappa$}
\label{Fig:LSR2_Elongation}
\end{subfigure}
\caption{Elongation of shortest curvature-bounded path $\gamma_m\in \Gamma(X,Y)$ with $(X,Y)\in \mathcal{O}_1\cup \mathcal{O}_2\cup \mathcal{O}_3$.}
\label{Fig:LSR_Elongation0}
\end{figure}

An elongation strategy for $(X,Y)\in \mathcal{O}_4 \cup \mathcal{O}_5$ has been provided by Fig.~\ref{Fig:d4_Elongation}. Regarding the case of $(X,Y)\in \mathcal{O}_1 \cup \mathcal{O}_2\cup \mathcal{O}_3$, an elongation strategy is  illustrated in Fig.~\ref{Fig:LSR_Elongation0}.

\begin{figure}[!htp]
\centering
\begin{subfigure}{4cm}
\centering
\includegraphics[width=0.7\linewidth]{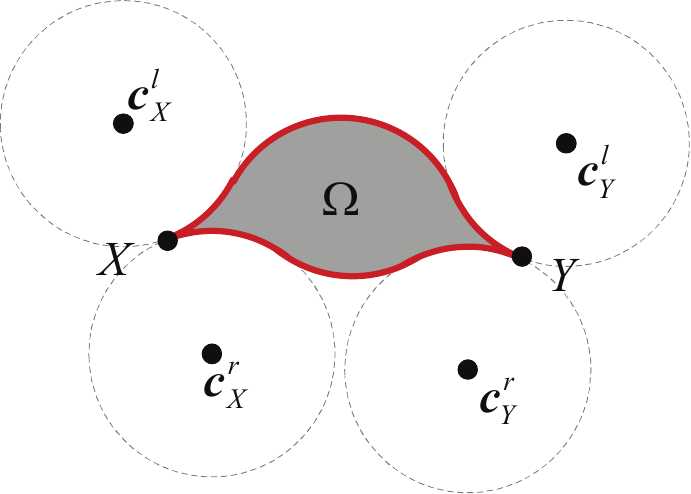}
\caption{$\gamma_m(X,Y)$ is of  RSR}
\label{Fig:Region_Omega_RSR}
\end{subfigure}
\begin{subfigure}{4cm}
\centering
\includegraphics[width=0.58\linewidth]{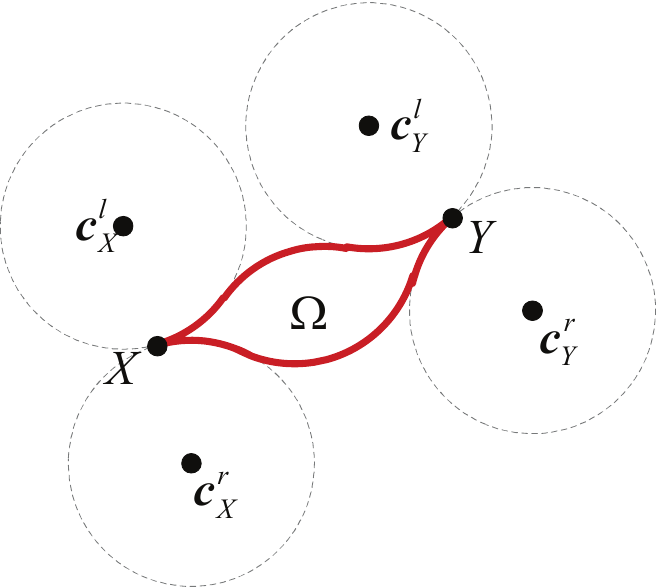}
\caption{$\gamma_m(X,Y)$ is of LSL}
\label{Fig:Region_Omega_LSL}
\end{subfigure}

\begin{subfigure}{4cm}
\centering
\includegraphics[width=0.7\linewidth]{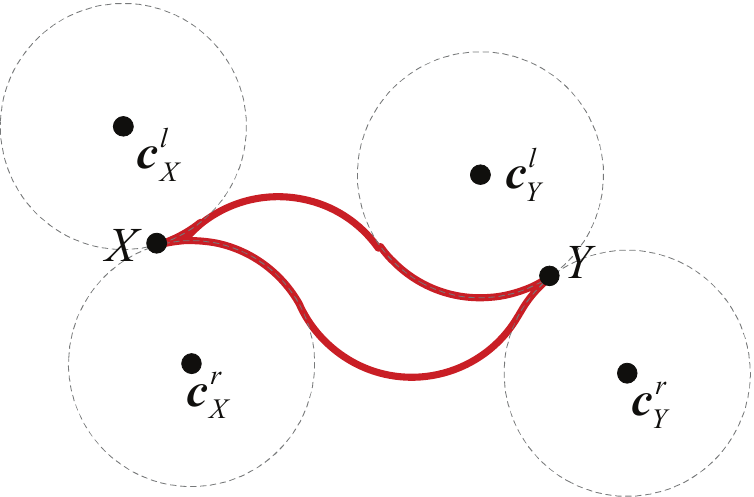}
\caption{$\gamma_m(X,Y)$ is of RSL}
\label{Fig:Region_Omega_RSL}
\end{subfigure}
\begin{subfigure}{4cm}
\centering
\includegraphics[width=0.4\linewidth]{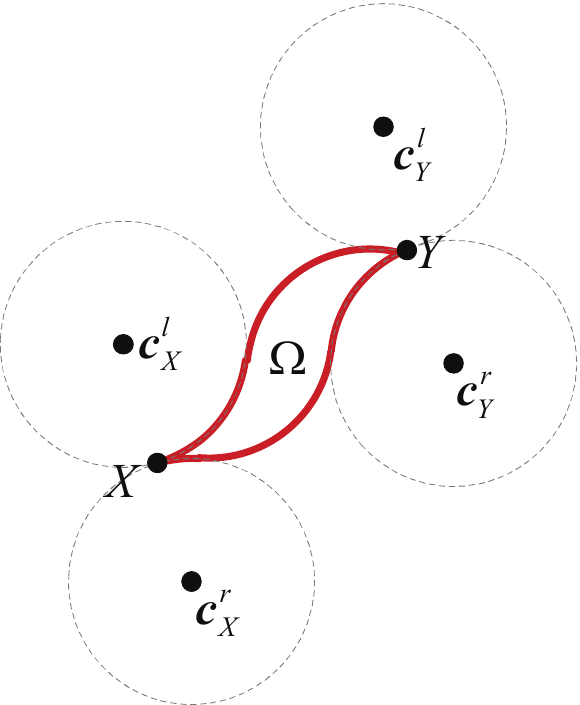}
\caption{$\gamma_m(X,Y)$ is of LSR}
\label{Fig:Region_Omega_LSR}
\end{subfigure}
\caption{Geometry for the closed region $\Omega$.}\label{Fig:Region_Omega}
\end{figure}

\begin{figure}[!htp]
\centering
\includegraphics[width = 4.0cm]{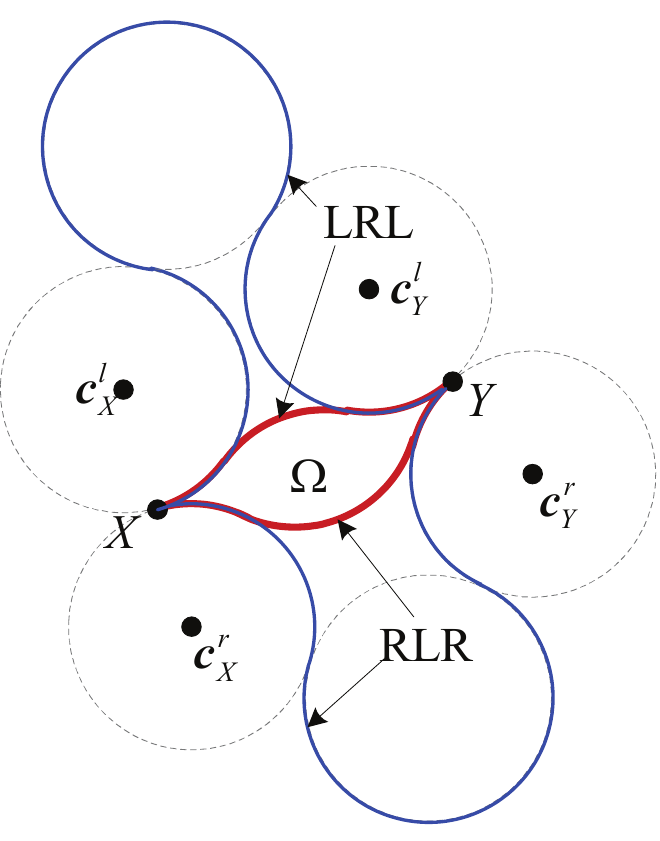}
\caption{Existene of two LRL-paths and two RLR-paths for $(X,Y)\in \nabla \mathcal{O}$.}
\label{Fig:Two_CCC}
\end{figure}

It has been shown in \cite{Ayala:2014,Ayala:2016} that  if $(X,Y)\in \nabla{\mathcal{O}}$ there exists a closed region $\Omega$, as shown in Fig.~\ref{Fig:Region_Omega}. Once $(X,Y)\in \nabla \mathcal{O}$, there exist two paths of type RLR and two paths of type LRL, as shown in Fig.~\ref{Fig:Two_CCC}. We denote by $\mathrm{RLR}^s$ and $\mathrm{RLR}^l$ the shorter and longer RLR-paths, respectively, and the same applies to $\mathrm{LRL}^s$ and $\mathrm{LRL}^l$.  It is apparent that the closed region $\Omega$ is bounded by the $\mathrm{RLR}^s$- and $\mathrm{LRL}^s$-paths.

Denote by $\ell_{\mathrm{LRL}}^s$, $\ell_{\mathrm{LRL}}^{l}$, $\ell_{\mathrm{RLR}}^s$, and $\ell_{\mathrm{RLR}}^{l}$ the lengths of  $\mathrm{RLR}^s$-, $\mathrm{RLR}^l$-, $\mathrm{LRL}^s$-, and $\mathrm{LRL}^l$-paths from $X$ to $Y$, respectively. Accordingly, we denote by $\ell_{\mathrm{RSR}}$, $\ell_{\mathrm{RSL}}$, $\ell_{\mathrm{LSR}}$, and $\ell_{\mathrm{LSL}}$ the lengths of   RSR-, RSL-, LSR-, and LSL-paths from $X$ to $Y$, respectively. Once RSR-path does not exist, we set $\ell_{\mathrm{RSR}} = +\infty$, and the same  applies to $\ell_{\mathrm{RSL}}$, $\ell_{\mathrm{LSR}}$,  $\ell_{\mathrm{LSL}}$, $\ell_{\mathrm{LRL}}^s$, $\ell_{\mathrm{LRL}}^{l}$, $\ell_{\mathrm{RLR}}^s$, and $\ell_{\mathrm{RLR}}^{l}$, accordingly.  
When $(X,Y)\in \nabla \mathcal{O}$ so that the closed region $\Omega$ exists, we set
\begin{align}
\begin{split}
\ell_1 \coloneqq &\  \max \{\ell_{\mathrm{LRL}}^s,\ \ell^s_{\mathrm{RLR}}\}\\
\ell_2  \coloneqq &\  \min 
\left\{ 
\begin{array}{l}
\ell_m+2\pi,\ \ell^l_{\mathrm{LRL}},\ \ell^l_{\mathrm{RLR}},\{\ell_{\mathrm{RSR}},\\
  \ell_{\mathrm{RSL}}, \ell_{\mathrm{LSR}}, \ell_{\mathrm{LSL}} \} \setminus \{\ell(\gamma_m)\}
\end{array}
\right\}
\end{split}
\label{EQ:S1_S2}
\end{align}
Given any $(X,Y)\in \nabla \mathcal{O}$, the lengths of all the CSC- and CCC-paths, once they exist, can be computed by geometric analysis \cite{Shkel:2001,Dubins:1957}. Thus,  both $\ell_1$ and $\ell_2$ in Eq.~(\ref{EQ:S1_S2}) can be readily obtained. According to \cite[{Proposition 1}]{Yao:2017}, we have $\ell_m < \ell_1$.  By the following lemma, we shall show that $\ell_1< \ell_2$ when the closed region $\Omega$ exists. 
\begin{lemma}
Given  any two oriented points $X$ and $Y$ so that $(X,Y) \in \nabla{\mathcal{O}}$, we have $\ell_1 < \ell_2$. 
\end{lemma}
Proof. By definition, for every C$_{\eta}$C$_{\xi}$C$_{\zeta}$-paths related to $\ell^s_{\mathrm{RLR}}$ and $\ell^s_{\mathrm{LRL}}$, the middle circular arc  is a minor arc, i.e, $\xi<\pi$. According to \cite[Lemma 2]{patsko2003three}, once the closed region $\Omega$ exists, the sum of the other two circular arcs is less than the middle circular arc, i.e., $\eta + \zeta < \xi$. Thus, we have $\ell_1 = \eta + \xi  + \zeta < 2\pi$, indicating $\ell_1< \ell_m + 2\pi$. If $\ell_2 = \ell^l_{\mathrm{LRL}}$ or $\ell^l_{\mathrm{RLR}}$, we immediately have $\ell_1< \ell_2$ by the definitions of $\ell^s_{\mathrm{RLR}}$, $\ell^s_{\mathrm{LRL}}$, $\ell^l_{\mathrm{RLR}}$, and $\ell^l_{\mathrm{LRL}}$. 

From now on, we compare $\ell_1$ with the element in  $\{\ell_{\mathrm{RSR}}, \ell_{\mathrm{RSL}}, \ell_{\mathrm{LSR}}, \ell_{\mathrm{LSL}} \} \setminus \{\ell(\gamma_m)\}$. Without loss of generality, let us assume that $\ell_{\mathrm{RSR}}\neq \ell_m$, i.e., $\ell_{\mathrm{RSR}}$ is an element of $\{\ell_{\mathrm{RSR}}, \ell_{\mathrm{RSL}}, \ell_{\mathrm{LSR}}, \ell_{\mathrm{LSL}} \} \setminus \{\ell(\gamma_m)\}$. Then, we have that the RSR-path has intersections with the boundary of $\Omega$. As a result, we have that $\ell_{\mathrm{RSR}}$ is greater than $\min\{\ell^l_{\mathrm{RLR}}, \mathrm{R}_{\mathrm{LRL}}^l\}$ according to \cite[Lemma 10]{Yao:2017}. The same applies to other CSC-paths with its length in $\{\ell_{\mathrm{RSR}}, \ell_{\mathrm{RSL}}, \ell_{\mathrm{LSR}}, \ell_{\mathrm{LSL}} \} \setminus \{\ell(\gamma_m)\}$, completing the proof.
 $\Box$
\begin{figure}
\centering
\begin{subfigure}{5cm}
\includegraphics[width = 5cm]{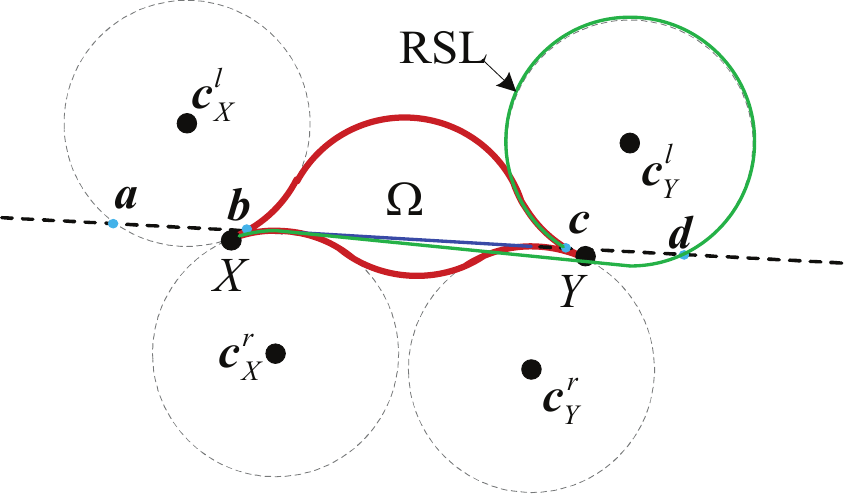}
\caption{$\mathrm{R_{\eta}S_d L_{\xi}}$ with the existence of $\Omega$}
\label{Fig:S_RSL}
\end{subfigure}
\begin{subfigure}{5cm}
\includegraphics[width = 5cm]{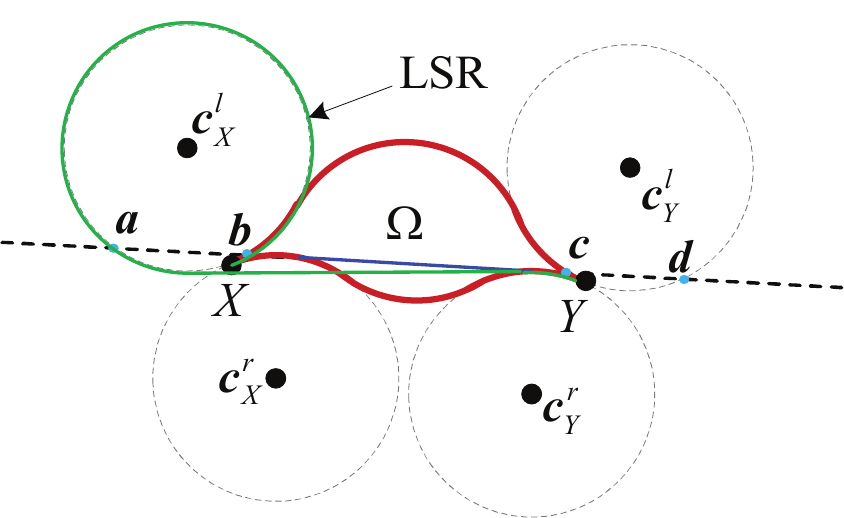}
\caption{$\mathrm{L_{\eta}S_d R_{\xi}}$ with the existence of $\Omega$}
\label{Fig:S_LSR}
\end{subfigure}
\begin{subfigure}{5cm}
\includegraphics[width = 5cm]{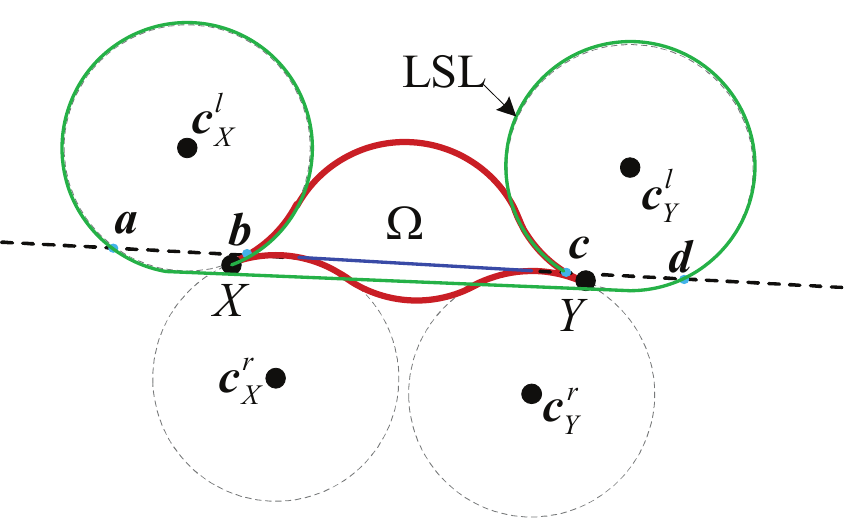}
\caption{$\mathrm{L_{\eta}S_d L_{\xi}}$ with the existence of $\Omega$}
\label{Fig:S_LSL}
\end{subfigure}
\caption{The CSC-paths with the existence of $\Omega$.}
\label{Fig:S_CSC}
\end{figure}
\begin{lemma}\label{LE:S_CSC}
Given any two oriented points $X$ and $Y$ so that $(X,Y)\in \nabla{\mathcal{O}}$, each $\mathrm{CSC}$-path, with its length in $\{\ell_{\mathrm{RSR}}, \ \ell_{\mathrm{RSL}},\ \ell_{\mathrm{LSR}},\ \ell_{\mathrm{LSL}}\}\setminus \{\ell_m\}$, has parallel tangents.
\end{lemma}
Proof. Without loss of generality, let us assume that the shortest path related to $\ell_m$ is of type RSR. Take the picture in Fig.~\ref{Fig:S_CSC} as an example. The extension of the straight line segment of the RSR-path has four points intersecting with the circles $C_X^l$ and $C_Y^l$, as shown by $\boldsymbol{a}$, $\boldsymbol{b}$, $\boldsymbol{c}$, and $\boldsymbol{d}$ in Fig.~\ref{Fig:S_CSC}. Because of the existence of the closed region $\Omega$, we have that the arcs $\wideparen{\boldsymbol{a}\boldsymbol{x}\boldsymbol{b}}$ and $\wideparen{\boldsymbol{c}\boldsymbol{y}\boldsymbol{d}}$ are minor arcs.
 Regarding the RSL-path, we have that the switching point from S to L lies on the arc $\wideparen{\boldsymbol{c}\boldsymbol{y}\boldsymbol{d}}$ (see Fig.~\ref{Fig:S_RSL}), indicating that the left-turn arc along the RSL-path  is a major arc. As for the LSR-path, we have that the switching point from L to S lies on the arc $\wideparen{\boldsymbol{a}\boldsymbol{x}\boldsymbol{b}}$ (see Fig.~\ref{Fig:S_LSR}), indicating that  the left-turn arc along the LSR-path is a major arc. As for the LSL-path, we have that the switching point from L to S and the switching point from S to L lie on $\wideparen{\boldsymbol{c}\boldsymbol{y}\boldsymbol{d}}$ and $\wideparen{\boldsymbol{a}\boldsymbol{x}\boldsymbol{b}}$, respectively, as shown by Fig.~\ref{Fig:S_LSL}, indicating that the both left-turn arcs are major arcs. Therefore, if the RSR-path is the shortest path, we have that at least one circular arc on each of the RSL-, LSR-, and LSL-paths  is a major arc, implying that all the RSL-, LSR-, and LSL-paths have parallel tangents. 

\begin{figure}[!htp]
\centering
\includegraphics[width = 4cm]{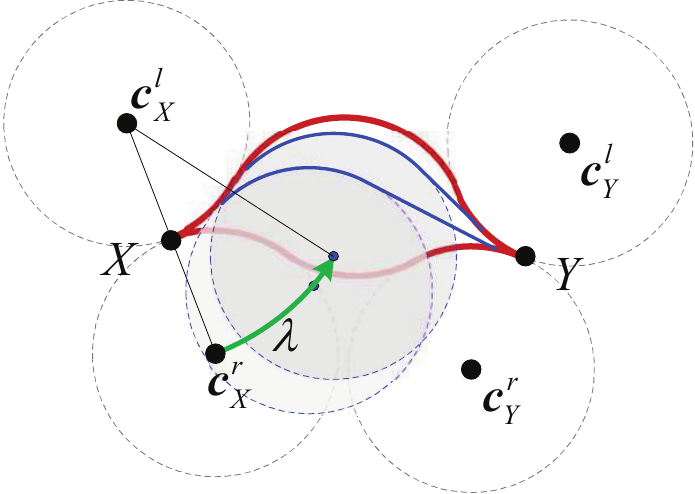}
\caption{Elongation of shortest curvature-bounded path $\gamma_m$ for $(X,Y)\in \nabla \mathcal{O}$.}
\label{Fig:RSR_Elongation_L1}
\end{figure}

Analogously, we can prove by the same way that if one of the RSL-, LSR-, and LSL-paths is the shortest path, then each other CSC-path with its length in  $\{\ell_{\mathrm{RSR}}, \ \ell_{\mathrm{RSL}},\ \ell_{\mathrm{LSR}},\ \ell_{\mathrm{LSL}}\}\setminus \{\ell_m\}$ has parallel tangents, completing the proof.
$\Box$
\begin{proposition}\label{LE:Lm_L1}
Given any two oriented points $X$ and $Y$ so that $(X,Y)\in \nabla{\mathcal{O}}$,  for any $s\in [\ell_m,\ell_1]\cup [\ell_2,+\infty)$ there exists a curvature-bounded path $\gamma  \in \Gamma(X,Y)$ so that $\ell(\gamma)  = s$.
\end{proposition}
Proof. We first consider the case that $s\in [\ell_m, \ell_1]$.  Let us consider to move a circular disk of radius $1/\kappa$ to deform the shortest curvature-bounded path $\gamma_m$, as shown in Fig.~\ref{Fig:RSR_Elongation_L1}.  By changing the value of $\lambda$, we can see that the shortest curvature-bounded path $\gamma_m$ can be elongated continuously from $\ell_m$ to $\ell_1$. Thus, for any $s\in [\ell_m, \ell_1]$ there exists a curvature-bounded path $\gamma \in \Gamma(X,Y)$ so that $s=\ell(\gamma)$. 

From now on, we consider the case of $s\in [\ell_2,+\infty)$. If $\ell_2 \in \{\ell^l_{\mathrm{RLR}}, \ell^l_{\mathrm{LRL}}\}$, then we have that the CCC-path with length of $\ell_2$ has parallel tangents according to the definitions of $\mathrm{RLR}^l$ and $\mathrm{LRL}^l$ in Fig.~\ref{Fig:Two_CCC}. Thus, if $\ell_2 \in \{ \ell^l_{\mathrm{RLR}}, \ell^l_{\mathrm{LRL}}\}$, the CCC-path with length of $\ell_2$ can be elongated to arbitrary length.  If $\ell_2 = \ell_m + 2\pi$, we have that the path has parallel tangents.  If $\ell_2\in \{\ell_{\mathrm{RSR}}, \ \ell_{\mathrm{RSL}},\ \ell_{\mathrm{LSR}},\ \ell_{\mathrm{LSL}}\}\setminus \{\ell_m\}$, we have that the path related to $\ell_2$ has parallel tangents, according to Lemma \ref{LE:S_CSC}. Therefore, the path related to $\ell_2$ can be elongated to arbitrary length, completing the proof.$\Box$

Proposition \ref{LE:Lm_L1} indicates that for any $s\in [\ell_m,\ell_1]\cup[\ell_2,+\infty)$ we can find a curvature-bounded path $\gamma \in \Gamma(X,Y)$ so that $s = \ell(\gamma)$. From now on, we will prove that for any  $s\in (\ell_1,\ell_2)$ it is impossible to find a curvature-bounded path $\gamma \in \Gamma(X,Y)$ so that $s = \ell(\gamma)$. 

Before proceeding, we present some useful notations.  Denote by $\mathcal{R}_{\boldsymbol{w}}(s)\subset \mathbb{R}^2$ the set that can be reached by all the curvature-bounded paths $\gamma$ of length $s$, starting from $X$ and ending with the final tangent vector being $\boldsymbol{w}$, i.e.,
\begin{align}
\mathcal{R}_{\boldsymbol{w}}(s) : = \{\boldsymbol{z}\in \mathbb{R}^2 | \ell(\gamma) = s \ \text{with}\ \gamma \in \Gamma(X,(\boldsymbol{z},\boldsymbol{w}))\}.\nonumber
\end{align}
As usual, we denote by $\partial \mathcal{R}_{\boldsymbol{w}}(s)$ the boundary of $\mathcal{R}_{\boldsymbol{w}}(s)$, and denote by $\mathrm{Int}\mathcal{R}_{\boldsymbol{w}}(s)$ the interior of $\mathcal{R}_{\boldsymbol{w}}(s)$. Denote by $P_{\mathrm{RSR}}^{\boldsymbol{w}}(s)\in \mathbb{R}^2$ the set of all the points that can be reached by RSR-path of length $s>0$, starting  from $X$ and ending with the final tangent being $\boldsymbol{w}$. The same explanation applies to $P_{\mathrm{RSL}}^{\boldsymbol{w}}(s)$, $P_{\mathrm{LSR}}^{\boldsymbol{w}}(s)$, $P_{\mathrm{LSL}}^{\boldsymbol{w}}(s)$, $P_{\mathrm{RLR}}^{\boldsymbol{w}}(s)$, and $P_{\mathrm{LRL}}^{\boldsymbol{w}}(s)$.
\begin{lemma}\label{LE:interior1}
Given any two oriented points $X$ and $Y$ so that $(X,Y) \in \nabla\mathcal{O}$, the following two statements hold:
\begin{description}
\item (1)   there exists a positive number $\delta > 0 $ so that  $\boldsymbol{y}\in \mathcal{R}_{\boldsymbol{w}}(\ell_1 - \varepsilon)$ for every $\varepsilon \in (0,\delta)$;
\item (2) $\boldsymbol{y}\not\in \mathcal{R}_{\boldsymbol{w}}(\ell_1 + \varepsilon)$  for any sufficiently small $\varepsilon > 0$.
\end{description}
\end{lemma}
Proof. By the definition of $\ell_1$ in Eq.~(\ref{EQ:S1_S2}), we have $\boldsymbol{y}\in \mathcal{R}_{\boldsymbol{w}}(\ell_1)$.  According to Proposition \ref{LE:Lm_L1}, for any $s\in [\ell_m,\ell_1]$, there exists a curvature-bounded path from $X$ to $Y$ with its   length being $s$. Thus, we have that $\boldsymbol{y}\in \mathcal{R}_{\boldsymbol{w}}(\ell_1-\eta)$ for $\eta \in [0,\ell_1-\ell_m]$. Since $\ell_1-\ell_m >0 $ (cf. \cite[{Proposition 1}]{Yao:2017}), it follows that  there exists $\delta>0$ so that $\boldsymbol{y}\in \mathcal{R}_{\boldsymbol{w}}(\ell_1-\varepsilon)$ for every $\varepsilon \in (0,\delta)$.

\begin{figure}[!htp]
\centering
\includegraphics[width = 10cm]{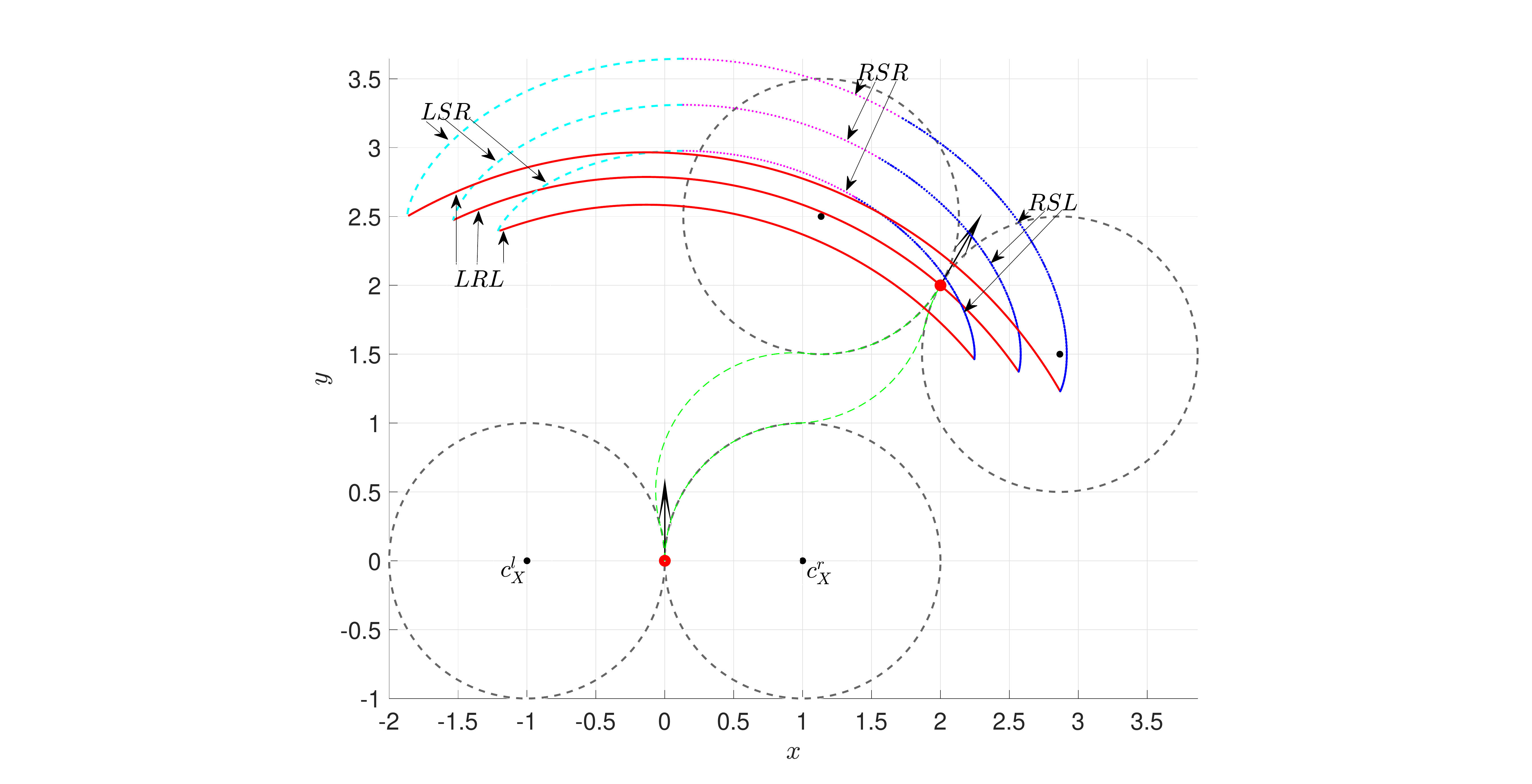}
\caption{Illustration of $\boldsymbol{y}\in P_{\mathrm{RLR}}^{\boldsymbol{w}}(\ell_1)\subset \mathcal{R}_{\boldsymbol{w}}(\ell_1)$.}
\label{Fig:AllPaths}
\end{figure}

From now on, we proceed to proving $\boldsymbol{y}\not \in \mathcal{R}_{\boldsymbol{w}}(\ell_1 + \varepsilon)$ for sufficiently small $\varepsilon$. 
Without loss of generality, let us assume that $\ell_1 = \ell_{\mathrm{RLR}}^s$, i.e., the CCC-path with length of $\ell_1$ from $X$ to $Y$ is of type RLR. In this case, according to  \cite{PatFed20}, we have that the set $P_{\mathrm{RLR}}^{\boldsymbol{w}}(\ell_1)$ belongs to the boundary of  $\mathcal{R}_{\boldsymbol{w}}(\ell_1)$, as shown by Fig.~\ref{Fig:AllPaths}. Notice that  $\boldsymbol{y}\in P_{\mathrm{RLR}}^{\boldsymbol{w}}(\ell_1)$, indicating $\boldsymbol{y}\in \partial \mathcal{R}_{\boldsymbol{w}}(\ell_1)$.

Let $\mathcal{N}_{\eta}(\boldsymbol{y})\subset \mathbb{R}^2$ be a circular neighborhood centered at $\boldsymbol{y}$ with radius $\eta > 0$, i.e.,
$$\mathcal{N}_{\eta}(\boldsymbol{y}) \coloneqq \{\boldsymbol{z}\in \mathbb{R}^2| \|\boldsymbol{z}-\boldsymbol{y}\|\leq \eta\}.$$
As $\boldsymbol{y}\in \partial \mathcal{R}_{\boldsymbol{w}}(\ell_1)$, we can choose a straight line segment $\boldsymbol{p}:[0,2\ell_1]\rightarrow \mathcal{N}_{\eta}(\boldsymbol{y})$   so that 
\begin{itemize}
\item $\boldsymbol{p}(\ell_1) = \boldsymbol{y}$,
\item  $\boldsymbol{p}(t)\not \in \mathcal{R}_{\boldsymbol{w}}(t)$ for $t\in [0,\ell_1)$, and 
\item $\boldsymbol{p}(t)  \in \mathcal{R}_{\boldsymbol{w}}(t)$ for $t\in (\ell_1,2\ell_1]$.
\end{itemize}
Let us consider a monotonically increasing sequence $(\varepsilon_i)_{i\in \mathbb{N}}$ so that $\varepsilon_0 = 0$ and $\varepsilon_i < 2\ell_1$ for $i\in \mathbb{N}$, and assume that there exists a large positive integer $N>0$ so that $\varepsilon_N = \ell_1$. If $\eta>0$ is small enough, for every $\varepsilon_i$, there exists a $\delta_i$ so that $\boldsymbol{p}(\varepsilon_i) \in P^{\boldsymbol{w}}_{\mathrm{RLR}}(\ell_1 + \delta_i)$.  Since the set $ P^{\boldsymbol{w}}_{\mathrm{RLR}}(\ell_1 + \delta_i)$ belongs to the boundary of $\mathcal{R}_{\boldsymbol{w}}(\ell_1 + \delta_i)$,  it follows $\boldsymbol{p}(\varepsilon_i)\in \partial \mathcal{R}_{\boldsymbol{w}}(\ell_1 + \delta_i)$. 

As the sequence $(\varepsilon_i)$ is monotonically increasing, the direction of the straight line $\boldsymbol{p}$ can be chosen so that the sequence $(\delta_i)$ is monotonically increasing as well.  Thus, picking an integer $\hat{i}\in \mathbb{N}$, we have $\boldsymbol{p}(\varepsilon_{\hat{i}+1}) \in \mathrm{Int} \mathcal{R}_{\boldsymbol{w}}(\ell_1 + \delta_{\hat{i}})$ according to the first statement of this lemma.  As $\boldsymbol{p}(\varepsilon_{\hat{i}}) \in \partial \mathcal{R}_{\boldsymbol{w}}(\ell_1+\delta_{\hat{i}})$, it follows $\boldsymbol{p}(\varepsilon_{\hat{i}-1}) \not \in  \mathcal{R}_{\boldsymbol{w}}(\ell_1 + \delta_{\hat{i}})$. If $\hat{i}=N+1$, we then have $\boldsymbol{p}(\varepsilon_N)\not \in \mathcal{R}_{\boldsymbol{w}}(\ell_1+\delta_{N+1})$.  Note that $\varepsilon_N = \ell_1$ and $\delta_N =0$, indicating $\boldsymbol{y} = \boldsymbol{p}(\varepsilon_N) \not \in \mathcal{R}_{\boldsymbol{w}}(\ell_1+\delta_{N+1})$.
If $\varepsilon_{N+1}$ is close enough to $\varepsilon_N$, there exists $\delta>0$ so that $\delta_{N+1}< \delta$, completing the proof.  
$\Box$
\begin{lemma}\label{LE:S2}
Given any two oriented points $X$ and $Y$ so that $(X,Y) \in \nabla \mathcal{O}$, there exists $\delta > 0$ so that $\boldsymbol{y}\not\in \mathcal{R}_{\boldsymbol{w}}(\ell_2-\varepsilon)$ and $\boldsymbol{y}\in \mathcal{R}_{\boldsymbol{w}}(\ell_2+\varepsilon)$ for every $\varepsilon \in (0,\delta)$.
\end{lemma}
Proof. By the definition of $\ell_2$ in Eq.~(\ref{EQ:S1_S2}), we have $\boldsymbol{y}\in \mathcal{R}_{\boldsymbol{w}}(\ell_2)$.  According to Proposition \ref{LE:Lm_L1}, for any $s\in [\ell_2,+\infty)$ there exists a curvature-bounded path from $X$ to $Y$ with its length being $s$. Thus, we have $\boldsymbol{y}\in \mathcal{R}_{\boldsymbol{w}}(\ell_2 + \eta)$ for any $\eta \geq 0$, indicating that there exists $\delta>0$ so that $\boldsymbol{y}\in \mathcal{R}_{\boldsymbol{w}}(\ell_2 + \varepsilon)$ for every $\varepsilon \in (0,\delta)$. 

From now on, we proceed to proving $\boldsymbol{y}\not\in \mathcal{R}_{\boldsymbol{w}}(\ell_2-\varepsilon)$ for sufficiently small $\varepsilon > 0$. By contradiction, let us assume that there exists a sufficiently small $\varepsilon >0$ so that $\boldsymbol{y}\in \mathcal{R}_{\boldsymbol{w}}(\ell_2 - \varepsilon)$.
 
Let us consider an object moving along a straight line $\boldsymbol{p}:[0,+\infty)$ so that at the instant $\ell_2$ the object reaches the point $\boldsymbol{y}$, i.e., $\boldsymbol{p}(\ell_2) = \boldsymbol{y}$. 
Let $\mathcal{N}_{\eta}(\boldsymbol{y})\subset \mathbb{R}^2$ be a circular neighborhood centered at $\boldsymbol{y}$ with radius $\eta > 0$, i.e.,
$$\mathcal{N}_{\eta}(\boldsymbol{y}) \coloneqq \{\boldsymbol{z}\in \mathbb{R}^2| \|\boldsymbol{z}-\boldsymbol{y}\|\leq \eta\}.$$
Without loss of generality, we consider $\ell_2 = \ell_{\mathrm{RSR}}$ and $\ell_m = \ell_{\mathrm{RSL}}$, i.e., the path associated with $\ell_2$ is of type RSR and the shortest path is of type RSL. Then, because the multi-valued set $P_{\mathrm{LSR}}(t)$ for $t>0$ is continuous (cf. \cite[Lemma 6]{Buzikov:2021}), it follows that for any sufficiently small $\varepsilon >0$ there exists $\eta>0$ so that the lenght of LSR-path from $X$ to a point in $\mathcal{N}_{\eta}(\boldsymbol{y})$ with the final tangent being $\boldsymbol{w}$ takes values in $(\ell_{\mathrm{LSR}}-\varepsilon,\ell_{\mathrm{LSR}}+\varepsilon)$. 
Without loss of generality, assume that the speed of the moving object is constant and small enough so that $\boldsymbol{p}(t)\in \mathcal{N}_{\eta}(\boldsymbol{y})$ for any $t\in [0,\ell_2]$. Then, for every $t\in [0,\ell_2]$, the LSR-path from $X$ to the point $\boldsymbol{p}(t)$ with the final tangent being $\boldsymbol{w}$ has a length in $(\ell_{\mathrm{LSR}}-\varepsilon,\ell_{\mathrm{LSR}}+\varepsilon)$. By the definition of $\ell_2$ in Eq.~(\ref{EQ:S1_S2}), we have $\ell_2 < \ell_{\mathrm{LSR}}$. Therefore, if $\varepsilon >0$ is small enough, the minimum time for a Dubins vehicle to intercept the moving object by following an LSR-path is greater than $\ell_2$. Analogously, we can prove that the minimum time for a Dubins vehicle to intercept the moving object by following LSL-path or one of CCC-paths  
is greater than $\ell_2$. Therefore, the minimum time for a Dubins vehicle from $X$ to intercept the moving object with the final tangent being $\boldsymbol{w}$ by following a CSC- or a CCC-path is $\ell_2$. This further indicates that the minimum time for a Dubins vehicle from $X$ to intercept the moving object with the final tangent being $\boldsymbol{w}$ by following a curvature-bounded path is $\ell_2$ (cf. \cite[Theorem 7]{Buzikov:2021}). However, the contradicting assumption indicates that  the minimum time for a Dubins vehicle from $X$ to intercept the moving object with the final tangent being $\boldsymbol{w}$ by following a curvature-bounded path is less than $\ell_2$. Hence, by contraposition, the proof is completed. 
$\Box$

As a result of Lemma \ref{LE:interior1} and Lemma \ref{LE:S2}, we immediately obtain the following result. 
\begin{proposition}\label{LE:non-existence}
Given any two oriented points $X$ and $Y$ so that $(X,Y) \in \nabla\mathcal{O}$, for every curvature-bounded path $\gamma \in \Gamma (X,Y)$ we have $\ell(\gamma)\not\in  (\ell_1,\ell_2)$.
\end{proposition}
Proof. By contradiction, let us assume that there exists a curvature-bounded path $\gamma \in \Gamma(X,Y)$ so that $\ell(\gamma)\in (\ell_1,\ell_2)$. Then, according to Lemmas \ref{LE:interior1} and \ref{LE:S2}, there are more than two homotopy classes, contradicting with the result in \cite{Ayala:2016} that there are two homotopy classes. Hence, by controposition, the proof is completed. 
$\Box$

This proposition indicates that if $(X,Y)\in \nabla \mathcal{O}$, the shortest curvature-bounded path $\gamma_m$ cannot be enlongated to arbitrary length. According to Propositions \ref{LE:CCC}--\ref{LE:non-existence}, we eventually have the main result summarized in the following theorem.
\begin{theorem}\label{TH:final}
Given any two oriented points $X$ and $Y$ so that $X\neq Y$, we have
\begin{description}
\item (1) If $(X,Y)\not \in \mathcal{O} \cup \nabla \mathcal{O}$, for  every $s\geq \ell_m$ there exists a curvature-bounded path $\gamma \in \Gamma(X,Y)$ so that $\ell(\gamma) = s$.
\item (2) If $(X,Y)\in \mathcal{O}$, for every $s\geq \ell_m$ there exists a curvature-bounded path $\gamma \in \Gamma(X,Y)$ so that $\ell(\gamma) = s$.
\item (3) If $(X,Y)\in \nabla\mathcal{O}$, the following two statesment hold:
\begin{itemize} 
\item for every $s\in [\ell_m,\ell_1]\cup [\ell_2,+\infty)$ there exists a curvature-bounded path $\gamma \in \Gamma(X,Y)$ so that $\ell(\gamma) = s$;
\item for  every $\gamma \in \Gamma(X,Y)$ we have $\ell(\gamma)\not \in (\ell_1,\ell_2)$. 
\end{itemize}
\end{description}
\end{theorem}
This theorem gives the necessary and sufficient conditions for the existence of curvature-bounded path with an expected length. Given any $X$ and $Y$ in $T\mathbb{R}^2$, the shortest curvature-bounded path $\gamma_m$ can be analytically obtained, indicating that the satisfications of $(X,Y)\in \mathcal{O}$ and $(X,Y)\in \nabla \mathcal{O}$ can be readily checked. In addition, the values of $\ell_1$ and $\ell_2$, once exist, can be computed analytically as well. Therefore, for any $X$ and $Y$ in $T\mathbb{R}^2$ all the conditions in Theorem \ref{TH:final} are numerically or analytically verifiable, allowing to predict the existance of curvature-bounded paths in $\Gamma(X,Y)$ with an expected length. 
\section{Numerical Examples}\label{SE:Numerical}
In this section, some examples of minimum-time path planning for  multiple fixed-wing  Unmanned Aerial Vehicles (UAVs) to simultaneously achieve a triangle-shaped flight formation will be presented to illustrate the developments of the paper. 

When considering that the fixed-wing UAVs fly in altitude hold mode with constant cruise speed,  the kinematics of such UAVs is the same as Dubins vehicle \cite{chen2020dubins}. Thus, in order to realize the minimum-time operation of achieving a desired flight formation for multiple UAVs, we should plan curvature-bounded path for each UAV so that all the UAVs  can reach their final conditions simultaneously. 

\begin{figure}[!htbp] 
\centering
\includegraphics[width = 4cm]{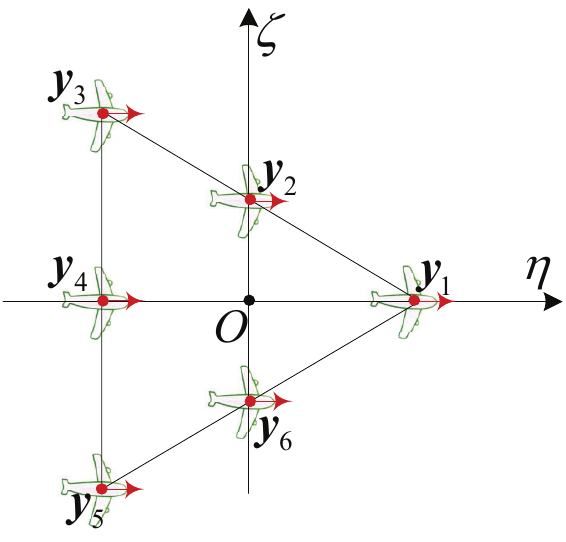}
\caption{Final positions of the six UAVs on the triangle-shaped flight formation.}
\label{Fig:Final_position_formation}
\end{figure}

We consider that there are $6$ UAVs, and the final position of UAV $\#i$ is denoted by $\boldsymbol{y}_i$ ($i=1,2,\ldots,6$).  Set the values of $\boldsymbol{y}_i$'s as
\begin{align}
\begin{tabular}{lll}
$\boldsymbol{y}_{1} = (\sqrt{3},0)$, &$\boldsymbol{y}_{2}=(0,1)$, & $\boldsymbol{y}_{3}=(-\sqrt{3},2)$\\
$\boldsymbol{y}_{4} = (-\sqrt{3},0)$, &$\boldsymbol{y}_{5}=(-\sqrt{3},-2)$, &$\boldsymbol{y}_{6} = (0,-1)$
\end{tabular}
\nonumber
\end{align}
Let the final tangent vectors of all the UAVs are the same, and are collinear with the $\eta$ axis of the $O\eta\zeta$ frame, as shown in Fig.~\ref{Fig:Final_position_formation}. It is clear from  Fig.~\ref{Fig:Final_position_formation} that  an equilateral-triangle formation will be formed if the six UAVs arrive their final conditions simultaneously.

Denote by $\boldsymbol{x}_i\in\mathbb{R}^2$ the initial position of UAV $\# i$, and let $\theta_i\in [0,2\pi)$ be the angle between the initial tangent vector and the $\eta$ axis, measured counterclockwise.  The values of $\boldsymbol{x}_i$'s and $\theta_i$'s are generated randomly by uniform distribution for three different cases, and are presented in Tables~\ref{tab1}--\ref{tab3}. 

\begin{table}[!htpb]
\centering
\caption{Case A: the values of $\boldsymbol{x}_i$'s and $\theta_i$'s.}
\label{tab1}
\begin{tabular}{ccc}
\toprule 
$i$ & $\boldsymbol{x}_i$ & $\theta_i$\\
   \midrule
$1$ & $(3.5313, -0.8619)$ & $0.5305$\\
$2$ & $(1.2238, 0.9698)$ & $4.8689$\\
$3$ & $(-3.5775, 1.3472)$ & $ 1.6328$\\
$4$ & $(1.6878, 0.9028)$ & $ 2.5119$\\
$5$ & $(2.9336, 0.0854)$ & $5.4582$\\
$6$ & $(1.1577, -0.0281)$ & $ 5.1353$\\
\bottomrule       
\end{tabular}
\end{table}
\begin{table}[!htpb]
\centering
\caption{Case B: the values of $\boldsymbol{x}_i$'s and $\theta_i$'s.}
\label{tab2}
\begin{tabular}{ccc}
\toprule 
$i$ & $\boldsymbol{x}_i$ & $\theta_i$\\
   \midrule
$1$ & $(4.3627, -1.0457)$ & $ 6.0141$\\
$2$ & $(-2.3376,0.2700)$ & $0.2919$\\
$3$ & $(2.1806, 3.3248)$ & $5.0283 $\\
$4$ & $(-0.8038, 3.4410)$ & $0.8915 $\\
$5$ & $(-4.5537, -1.3816)$ & $2.6500$\\
$6$ & $(1.5350, -0.2869)$ & $ 5.7537$\\
\bottomrule       
\end{tabular}
\end{table}
\begin{table}[!htpb]
\centering
\caption{Case C: the values of $\boldsymbol{x}_i$'s and $\theta_i$'s.}
\label{tab3}
\begin{tabular}{ccc}
\toprule 
$i$ & $\boldsymbol{x}_i$ & $\theta_i$\\
   \midrule
$1$ & $(1.8829,4.4956)$ & $0.7477$\\
$2$ & $(-0.9264,0.0596)$ & $3.1313$\\
$3$ & $(-3.1202,1.1104)$ & $6.0302$\\
$4$ & $(-4.4641,1.4021)$ & $0.5136$\\
$5$ & $(1.6253, -3.9714)$ & $3.9600$\\
$6$ & $(-0.7889,-2.7028)$ & $1.4063$\\
\bottomrule       
\end{tabular}
\end{table}
Note that the initial tangent vector for UAV $\# i$ is given by $\boldsymbol{v}_i \coloneqq (\cos\theta_i,\sin \theta_i)$, and the final tangent vectors for all the six UAVs is the same as $\boldsymbol{w} = (1,0)$.  We denote by $\ell_m^i$ the length of the shortest curvature-bounded path for UAV $\# i$ from the initial oriented point $(\boldsymbol{x}_i,\boldsymbol{v}_i)$ to the final oriented point $(\boldsymbol{y}_{i},\boldsymbol{w})$. Accordingly, we denote by $\ell_1^i$ and $\ell_2^i$ the lengths corresponding to $\ell_1$ and $\ell_2$ defined in Eq.~(\ref{EQ:S1_S2}) for UAV $\# i$, i.e.,  in the case of $((\boldsymbol{x}_i,\boldsymbol{v}_i),(\boldsymbol{y}_{i},\boldsymbol{w}))\in \nabla \mathcal{O}$, there exists a curvature-bounded path  from $(\boldsymbol{x}_i,\boldsymbol{v}_i)$ to $(\boldsymbol{y}_{i},\boldsymbol{w})) $ if and only if its length lies in  $[\ell_m^i, \ell_1^i] \cup [\ell_2^i,+\infty) $. The values of $\ell_m^i$'s, $\ell_1^i$'s, and $\ell_2^i$'s for cases A, B, and C are presented in Table~\ref{tab4}. 
\begin{table}[htbp] 
\caption{The values of  $\ell_m^i$'s, $\ell_1^i$'s, and $\ell_2^i$'s for  cases A, B, and C.}
\label{tab4}
\centering
    \begin{tabular}{cccc}
\toprule 
         Item  & Case A  & Case B& Case C  \\ 
\midrule
        $\ell^1_m$ & 7.3871 & 8.5854 & 8.0845 \\   
        $\ell^1_1$ & $+\infty$ & $+\infty$ & $+\infty$  \\
         $\ell^1_2$ & $+\infty$ & $+\infty$ & $+\infty$ \\
\hline
        $\ell^2_m$ & 5.7164 & 2.4540 & 5.9104 \\      
      $\ell^2_1$ & $+\infty$ & 2.7219 & $+\infty$  \\
        $\ell^2_2$ & $+\infty$ & 8.7279 & $+\infty$ \\
 \hline
        $\ell^3_m$ & 7.0162 & 8.6103 & 7.8796 \\      
        $\ell^3_1$ & $+\infty$ & $+\infty$ & $+\infty$  \\
        $\ell^3_2$ & $+\infty$ & $+\infty$ & $+\infty$ \\
\hline
        $\ell^4_m$ & 6.7435 & 8.4646 & 3.3402 \\      
       $\ell^4_1$ & $+\infty$ & $+\infty$ & 3.6783  \\
        $\ell^4_2$ & $+\infty$ & $+\infty$ & 7.8609 \\
\hline
       $\ell^5_m$ & 9.7219 & 6.3674 & 6.6030 \\      
       $\ell^5_1$ & $+\infty$ & $+\infty$ & $+\infty$  \\
        $\ell^5_2$ & $+\infty$ & $+\infty$ & $+\infty$ \\
\hline
         $\ell^6_m$ & 6.7160 & 7.0891 & 7.6161 \\      
       $\ell^6_1$ & $+\infty$ & $+\infty$ & $+\infty$  \\
        $\ell^6_2$ & $+\infty$ & $+\infty$ & $+\infty$ \\
\bottomrule  
    \end{tabular}
\end{table}

Set 
$$\Phi_i \coloneqq
\begin{cases}
[\ell_m^i,\ell_1^i]\cup [\ell_2^i,+\infty)  \ \ \ \text{if} \ \ell_1^i < +\infty\ \text{and}\ \ell_2^i < +\infty\\
[\ell_m^i,+\infty)  \ \ \ \ \ \ \ \ \ \ \ \ \ \ \ \  \text{if} \ \ell_1^i = +\infty\ \text{and}\ \ell_2^i = +\infty\\
\end{cases} $$
Then, according to Theorem \ref{TH:final}, we have that the minimum time for the 6 UAVs to realize the triangle-shaped flight formation in Fig.~\ref{Fig:Final_position_formation} is given by
\begin{align}
t_m \coloneqq \min \{\Phi_1\cap\Phi_2\cap\Phi_3\cap\Phi_4\cap\Phi_5\cap\Phi_6 \}
\label{EQ:tm}
\end{align}
Combining Eq.~(\ref{EQ:tm}) and the data in Table~\ref{tab4}, we have that the minimum time to realize the triangle-shaped formation is $9.7219$, $8.7279$, and $8.0845$ for cases A, B, and C, respectively. 

\begin{figure}[!htbp]
    \centering
    \begin{subfigure}[b]{0.48\textwidth}
           \centering
           \includegraphics[width=\textwidth]{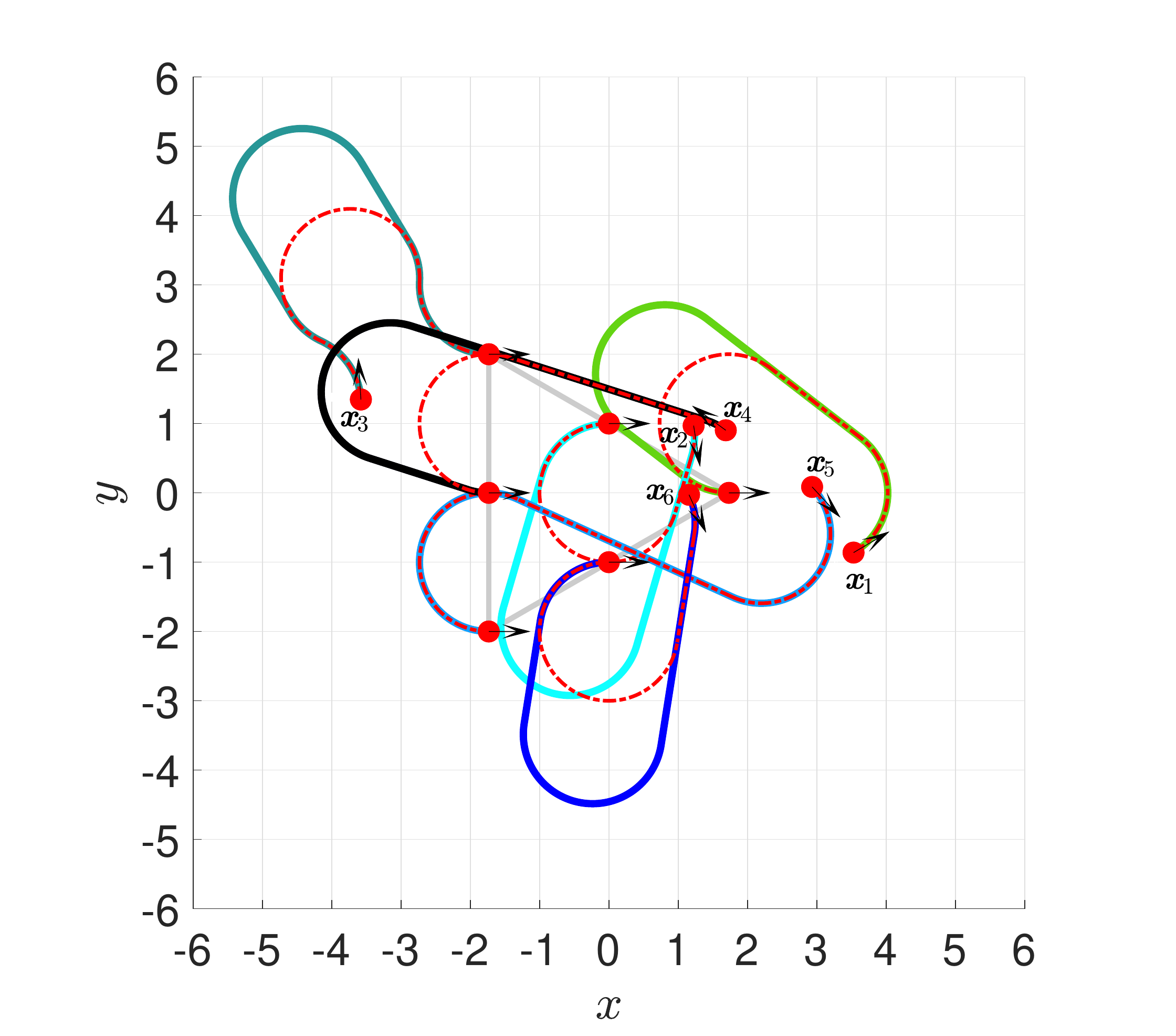}
            \caption{Case A}
\label{Fig:caseA_simulation}
    \end{subfigure}
    \begin{subfigure}[b]{0.48\textwidth}
            \centering
            \includegraphics[width=\textwidth]{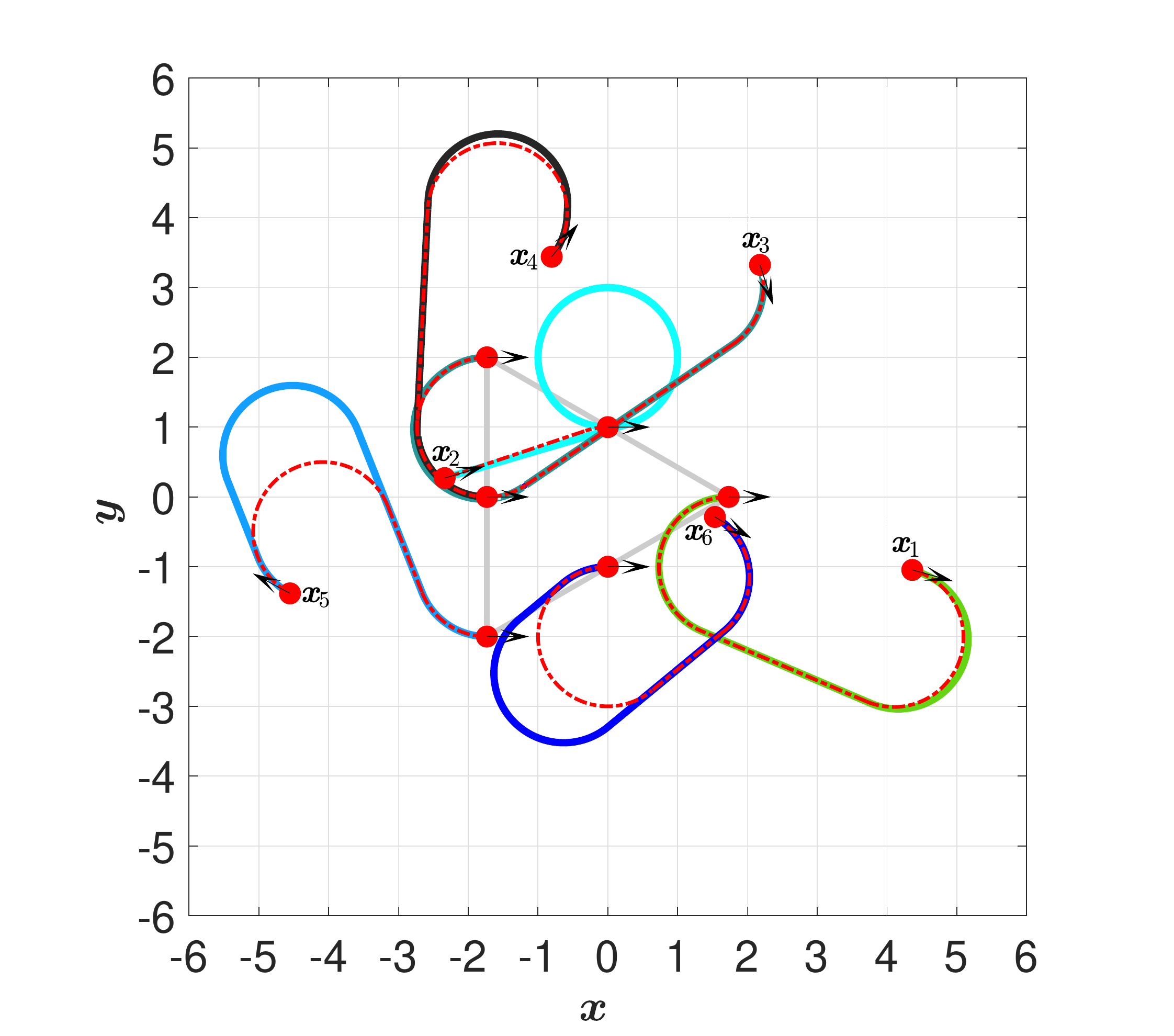}
\caption{Case B}\label{Fig:caseB_simulation}
    \end{subfigure}
    \begin{subfigure}[b]{0.48\textwidth}
            \centering
            \includegraphics[width=\textwidth]{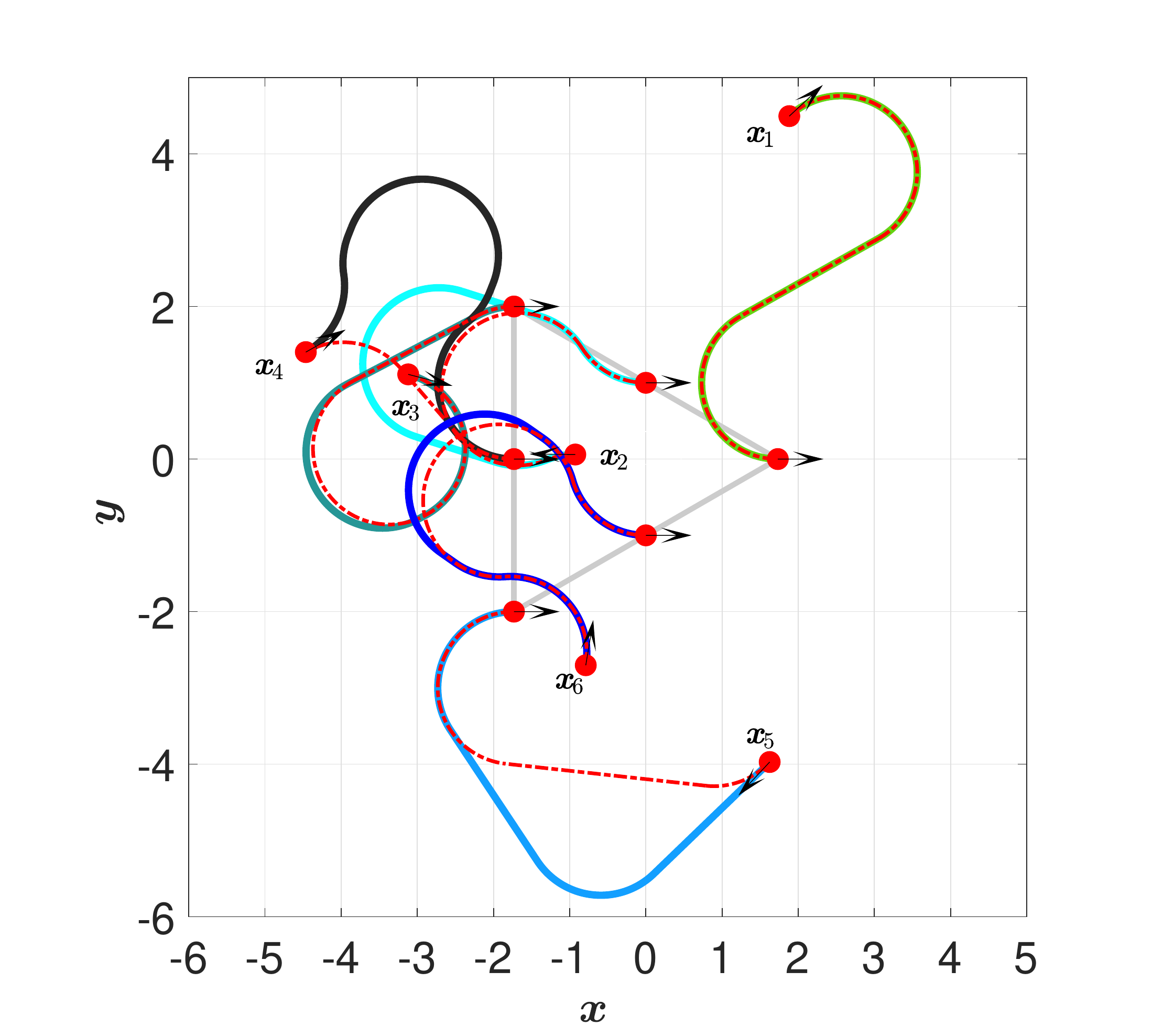}
\caption{Case C}\label{Fig:caseC_simulation}
    \end{subfigure}
    \caption{The cooperative paths for six UAVs to achieve a triangle-shaped flight formation.}
\label{Fig:simulation}
\end{figure}

Regarding case A, we can see from the data in Table~\ref{tab4} that in order to realize the minimum-time formation, we have to elongate the shortest paths for UAV $\#1$, UAV $\#2$, and UAV $\#3$, UAV $\#4$, and UAV $\#6$ to the length of $ 9.7219$.  The elongation strategies introduced in Figs.~\ref{Fig:CCC_Elongation}, \ref{Fig:d4_Elongation}, \ref{Fig:LSR_Elongation0}, and \ref{Fig:RSR_Elongation_L1} are used to elongate the shortest paths, accordingly. 
The elongated paths for UAV $\#1$, UAV $\#2$, and UAV $\#3$, UAV $\#4$, and UAV $\#6$  are presented by the solid curves  in Fig.~\ref{Fig:caseA_simulation}, where the dashed curves denote the shortest curvature-bounded paths.  If each UAV  follows its elongated path, all the six UAVs will achieve  a triangle-shaped flight formation in a minimum time. 

Regarding case B, we can see from the data in Table~\ref{tab4} that  in order to realize the minimum-time formation, we have to elongate the shortest paths for UAV $\#1$, UAV $\#3$,   UAV $\#4$, UAV $\#5$, and UAV $\#6$ to the length of $ 8.7279$. Note that UAV $\# 2$ follows a CSC-path with length of $\ell_2$ instead of the shortest path. The elongated path of each UAV for the minimum-time flight formation is presented by the solid curve in  Fig.~\ref{Fig:caseB_simulation}.

For case C, the paths of all the six UAVs have to be elongated to the length of $8.0845$ in order to realize the minimum-time operation of simultaneously achieving their final conditions, required for the triangle-shaped flight formation. Since $\ell_2^4 < 8.0845$, it follows that the shortest curvature-bounded path for UAV $\#4$ cannot be elongated to $8.0845$  according to Theorem \ref{TH:final}. We can get a curvature-bounded path of length by elongating  from the CCC-path of length $\ell_2^4 = 7.8609$, as shown by the solid curve starting from $\boldsymbol{x}_4$ in Fig.~\ref{Fig:caseC_simulation}. 
\section{Conclusions}\label{SE:Conclusions}
One of the fundamental issues in the field of cooperative guidance and path planning  for mulitple nonholonomic vehicles is to find curvature-bounded paths between two oriented points with  expected lengths. As the shortest curvature-bounded path between two oriented points can be analytically obtained, one usually elongates the shortest curvature-bounded path to an expected length. It was shown in the paper that, if shortest curvature-bounded path between two oriented points takes a geometric pattern of CCC, it can be elongated to arbitrary length. If the shortest curvature-bounded path between two oriented points takes a geometric pattern of CSC, the space of initial and final oriented points was divided into two subspaces $\mathcal{O}$ and $\nabla \mathcal{O}$. If the two oriented points lie in $\mathcal{O}$, the shortest curvature-bounded path can be elongated to arbitrary length. However, if  the two oriented points lie in $\nabla\mathcal{O}$, there exists an interval so that between the two oriented points  there is not a curvature-bounded path with its length  in the interval. 
In addition, the boundary values of the non-existence interval were presented. As a result, given any two oriented points and any expected length, it is readily to check if there exists a curvature-bounded path with the expected length; once it exists, an elongation strategy was provided to get the curvature-bounded path with the expected length. All the developments were finally demonstrated and verified by three examples of cooperative path planning for multiple UAVs to realize a triangle-shaped flight formation. 
\section*{Acknowledgements}
This work was supported by the National Natural Science Foundation of China (Grant Nos. 61903331 and 62088101).
\bibliographystyle{unsrt}  
\bibliography{references} 

\begin{thebibliography}{10}

\bibitem{markov1887some}
Andrey~Andreyevich Markov.
\newblock Some examples of the solution of a special kind of problem on
  greatest and least quantities.
\newblock {\em Soobshch. Karkovsk. Mat. Obshch}, 1:250--276, 1887.

\bibitem{Dubins:1957}
Lester~E. Dubins.
\newblock On curves of minimal length with a constraint on average curvature,
  and with prescribed initial and terminal positions and tangents.
\newblock {\em American Journal of Mathematics}, 79(3):497--516, 1957.

\bibitem{pontryagin1987mathematical}
Lev~Semenovich Pontryagin.
\newblock {\em Mathematical theory of optimal processes}.
\newblock CRC press, 1987.

\bibitem{sussmann1991shortest}
H{\'e}ctor~J Sussmann and Guoqing Tang.
\newblock Shortest paths for the {Reeds-Shepp} car: a worked out example of the
  use of geometric techniques in nonlinear optimal control.
\newblock {\em Rutgers Center for Systems and Control Technical Report},
  10:1--71, 1991.

\bibitem{boissonnat1994shortest}
Jean-Daniel Boissonnat, Andr{\'e} C{\'e}r{\'e}zo, and Juliette Leblond.
\newblock Shortest paths of bounded curvature in the plane.
\newblock {\em Journal of Intelligent and Robotic Systems}, 11(1):5--20, 1994.

\bibitem{Bui:1994}
Xuan-Nam Bui, Jean-Daniel {Boissonnat}, Philippe {Sou\`eres}, and Jean-Paul
  {Laumond}.
\newblock Shortest path synthesis for {Dubins} non-holonomic robot.
\newblock In {\em Proceedings of the 1994 IEEE International Conference on
  Robotics and Automation}, pages 2--7 vol.1, 1994.

\bibitem{boissonnat1996polynomial}
Jean-Daniel Boissonnat and Sylvain Lazard.
\newblock A polynomial-time algorithm for computing a shortest path of bounded
  curvature amidst moderate obstacles.
\newblock In {\em Proceedings of the twelfth annual symposium on Computational
  geometry}, pages 242--251, 1996.

\bibitem{Shkel:2001}
Andrei~M. Shkel and Vladimir Lumelsky.
\newblock Classification of the {Dubins} set.
\newblock {\em Robotics and Autonomous Systems}, 34(4):179 -- 202, 2001.

\bibitem{patsko2003three}
Valerii~S. Patsko, S.~G. Pyatko, and Andrey~A. Fedotov.
\newblock Three-dimensional reachability set for a nonlinear control system.
\newblock {\em Journal of Computer and Systems Sciences International},
  42(3):320--328, 2003.

\bibitem{schumacher2003path}
Corey Schumacher, Phillip Chandler, Steven Rasmussen, and David Walker.
\newblock Path elongation for {UAV} task assignment.
\newblock In {\em AIAA Guidance, Navigation, and Control Conference and
  Exhibit}, page 5585, 2003.

\bibitem{shanmugavel2005path}
Madhavan Shanmugavel, Antonios Tsourdos, Rafal Zbikowski, and Brian White.
\newblock Path planning of multiple {UAVs} using {Dubins} sets.
\newblock In {\em AIAA Guidance, Navigation, and Control Conference and
  Exhibit}, page 5827, 2005.

\bibitem{ma2006receding}
Xiang Ma and David~A Castanon.
\newblock Receding horizon planning for {Dubins} traveling salesman problems.
\newblock In {\em Proceedings of the 45th IEEE Conference on Decision and
  Control}, pages 5453--5458. IEEE, 2006.

\bibitem{mcgee2007optimal}
Timothy~G McGee and J~Karl Hedrick.
\newblock Optimal path planning with a kinematic airplane model.
\newblock {\em Journal of guidance, control, and dynamics}, 30(2):629--633,
  2007.

\bibitem{giordano2009shortest}
Paolo~Robuffo Giordano and Marilena Vendittelli.
\newblock Shortest paths to obstacles for a polygonal {Dubins} car.
\newblock {\em IEEE Transactions on Robotics}, 25(5):1184--1191, 2009.

\bibitem{techy2009minimum}
Laszlo Techy and Craig~A Woolsey.
\newblock Minimum-time path planning for unmanned aerial vehicles in steady
  uniform winds.
\newblock {\em Journal of guidance, control, and dynamics}, 32(6):1736--1746,
  2009.

\bibitem{SHANMUGAVEL20101084}
Madhavan Shanmugavel, Antonios Tsourdos, Brian White, and Rafał Żbikowski.
\newblock Co-operative path planning of multiple {UAVs} using {Dubins} paths
  with clothoid arcs.
\newblock {\em Control Engineering Practice}, 18(9):1084--1092, 2010.

\bibitem{ortiz2013multi}
Andres Ortiz, Derek Kingston, and C{\'e}dric Langbort.
\newblock Multi-{UAV} velocity and trajectory scheduling strategies for target
  classification by a single human operator.
\newblock {\em Journal of Intelligent \& Robotic Systems}, 70(1):255--274,
  2013.

\bibitem{savkin_hoy_2013}
Andrey~V. Savkin and Michael Hoy.
\newblock Reactive and the shortest path navigation of a wheeled mobile robot
  in cluttered environments.
\newblock {\em Robotica}, 31(2):323–330, 2013.

\bibitem{bakolas2013optimal}
Efstathios Bakolas and Panagiotis Tsiotras.
\newblock Optimal synthesis of the {Zermelo}--{Markov}--{Dubins} problem in a
  constant drift field.
\newblock {\em Journal of Optimization Theory and Applications},
  156(2):469--492, 2013.

\bibitem{Ayala:2014}
José Ayala and Hyam Rubinstein.
\newblock Non-uniqueness of the homotopy class of bounded curvature paths.
\newblock {\em arXiv}, 2014.

\bibitem{vavna2015dubins}
Petr V{\'a}{\v{n}}a and Jan Faigl.
\newblock On the {Dubins} traveling salesman problem with neighborhoods.
\newblock In {\em 2015 IEEE/RSJ International Conference on Intelligent Robots
  and Systems (IROS)}, pages 4029--4034. IEEE, 2015.

\bibitem{meyer2015dubins}
Yizhaq Meyer, Pantelis Isaiah, and Tal Shima.
\newblock On {Dubins} paths to intercept a moving target.
\newblock {\em Automatica}, 53:256--263, 2015.

\bibitem{Ayala:2015}
José Ayala.
\newblock Length minimising bounded curvature paths in homotopy classes.
\newblock {\em Topology and its Applications}, 193:140 -- 151, 2015.

\bibitem{Ayala:2016}
José Ayala and Hyam Rubinstein.
\newblock The classification of homotopy classes of bounded curvature paths.
\newblock {\em Israel Journal of Mathematics}, 213:79--107, 2016.

\bibitem{kaya2017markov}
C~Yal{\c{c}}{\i}n Kaya.
\newblock {Markov}--{Dubins} path via optimal control theory.
\newblock {\em Computational Optimization and Applications}, 68(3):719--747,
  2017.

\bibitem{Ayala+2017+283+292}
José Ayala.
\newblock On the topology of the spaces of curvature constrained plane curves.
\newblock {\em Advances in Geometry}, 17(3):283--292, 2017.

\bibitem{Yao:2017}
Weiran Yao, Naiming Qi, Jun Zhao, and Neng Wan.
\newblock Bounded curvature path planning with expected length for {Dubins}
  vehicle entering target manifold.
\newblock {\em Robotics and Autonomous Systems}, 97:217 -- 229, 2017.

\bibitem{vavna2018dubins}
Petr V{\'a}{\v{n}}a, Jakub Sl{\'a}ma, and Jan Faigl.
\newblock The {Dubins} traveling salesman problem with neighborhoods in the
  three-dimensional space.
\newblock In {\em 2018 IEEE International Conference on Robotics and Automation
  (ICRA)}, pages 374--379. IEEE, 2018.

\bibitem{manyam2019optimal}
Satyanarayana~Gupta Manyam, David Casbeer, Alexander Von~Moll, and Zachariah
  Fuchs.
\newblock Optimal {Dubins} paths to intercept a moving target on a circle.
\newblock In {\em 2019 American Control Conference (ACC)}, pages 828--834.
  IEEE, 2019.

\bibitem{chen2019shortest}
Zheng Chen and Tal Shima.
\newblock Shortest {Dubins} paths through three points.
\newblock {\em Automatica}, 105:368--375, 2019.

\bibitem{ding2019curvature}
Yulong Ding, Bin Xin, and Jie Chen.
\newblock Curvature-constrained path elongation with expected length for
  {Dubins} vehicle.
\newblock {\em Automatica}, 108:108495, 2019.

\bibitem{chen2019relaxed}
Zheng Chen and Tal Shima.
\newblock Relaxed {Dubins} problems through three points.
\newblock In {\em 2019 27th Mediterranean Conference on Control and Automation
  (MED)}, pages 501--506. IEEE, 2019.

\bibitem{Yao:2020}
Weiran {Yao}, Naiming {Qi}, Chengfei {Yue}, and Neng {Wan}.
\newblock Curvature-bounded lengthening and shortening for restricted vehicle
  path planning.
\newblock {\em IEEE Transactions on Automation Science and Engineering},
  17(1):15--28, 2020.

\bibitem{Yao:2020Trajectory}
Weiran Yao, Liming Xin, Yan Peng, Yang Chen, Naiming Qi, and Yu~Sun.
\newblock Trajectory consensus for coordination of multiple curvature-bounded
  vehicles.
\newblock {\em IEEE Transactions on Cybernetics}, pages 1--13, 2020.

\bibitem{PatFed20}
Valerii~S. Patsko and Andrey~A. Fedotov.
\newblock Analytic description of a reachable set for the {Dubins} car.
\newblock {\em Trudy Instituta Matematiki i Mekhaniki URO RAN}, 26(1):182--197,
  2020.

\bibitem{chen2020dubins}
Zheng Chen.
\newblock On {Dubins} paths to a circle.
\newblock {\em Automatica}, 117:108996, 2020.

\bibitem{jha2020shortest}
Bhargav Jha, Zheng Chen, and Tal Shima.
\newblock On shortest {Dubins} path via a circular boundary.
\newblock {\em Automatica}, 121:109192, 2020.

\bibitem{MATVEEV2020108831}
Alexey~S. Matveev, Valentin~V. Magerkin, and Andrey~V. Savkin.
\newblock A method of reactive control for 3d navigation of a nonholonomic
  robot in tunnel-like environments.
\newblock {\em Automatica}, 114:108831, 2020.

\bibitem{Buzikov:2021}
Maksim Buzikov and Andrey Galyaev.
\newblock Minimum-time lateral interception of a moving target by a {Dubins}
  car.
\newblock {\em Arxiv}, pages 1--16, 2021.

\bibitem{zheng2021time}
Yuan Zheng, Zheng Chen, Xueming Shao, and Wenjie Zhao.
\newblock Time-optimal guidance for intercepting moving targets by {Dubins}
  vehicles.
\newblock {\em Automatica}, 128:109557, 2021.

\bibitem{zheng2021time_angle}
Yuan Zheng and Zheng Chen.
\newblock Time-optimal guidance for intercepting moving targets with
  impact-angle constraints.
\newblock {\em Chinese Journal of Aeronautics (in Press)}, 2021.

\bibitem{chen2021descent}
Zheng Chen, Chen-hao Sun, Xue-ming Shao, and Wen-jie Zhao.
\newblock A descent method for the {Dubins} traveling salesman problem with
  neighborhoods.
\newblock {\em Frontiers of Information Technology \& Electronic Engineering},
  22(5):732--740, 2021.

\end{thebibliography}

\end{document}